\documentclass[12pt,a4j, twoside]{amsart}

\usepackage{docmute}
\usepackage{amsmath,amssymb}
\usepackage{amsthm}
\usepackage[abbrev]{amsrefs}
\usepackage{latexsym}
\usepackage{mathrsfs}

\theoremstyle{plain}
\newtheorem{thm}{Theorem}[section]
\newtheorem{prop}[thm]{Proposition}
\newtheorem{lem}[thm]{Lemma}
\newtheorem{cor}[thm]{Corollary}
\newtheorem{claim}[thm]{Claim}

\theoremstyle{definition}
\newtheorem{dfn}[thm]{Definition}
\newtheorem{ex}[thm]{Example}

\theoremstyle{remark}
\newtheorem{rem}[thm]{Remark}

\newtheorem*{ack}{Acknowledgement}

\renewcommand{\labelenumi}{(\theenumi)}
\numberwithin{equation}{section}

\newcommand{\R}{\mathbb{R}}
\newcommand{\N}{\mathbb{N}}

\newcommand{\F}{\mathscr{F}}

\newcommand{\supp}{\mathop{\mathrm{supp}}}
\newcommand{\diam}{\mathop{\mathrm{diam}}}
\newcommand{\midd}{\mathrel{} \middle| \mathrel{}}

\newcommand{\kf}{d_{\mathrm{KF}}}
\newcommand{\prok}{d_{\mathrm{P}}}

\newcommand{\lprok}[1]{d_{\mathrm{P}}^{(#1)}}
\newcommand{\Lip}{{\mathcal{L}}ip}
\newcommand{\conc}{d_{\mathrm{conc}}}
\newcommand{\haus}{d_{\mathrm{H}}}
\newcommand{\id}{\mathrm{id}}
\newcommand{\pr}{\mathrm{pr}}

\DeclareMathOperator{\OD}{ObsDiam}
\DeclareMathOperator{\PD}{diam}
\DeclareMathOperator{\LR}{LeRad}
\DeclareMathOperator{\df}{def}
\DeclareMathOperator{\lm}{lm}
\DeclareMathOperator{\arcsinh}{arcsinh}

\newcommand{\fracinline}[2]{
\raisebox{1.0ex}{ $#1$}
\raisebox{0ex}{\Large $/$}
\raisebox{-0.5ex}{\small $#2$}
}

\title{Concentration of product spaces}
\author{Daisuke Kazukawa}
\address{Mathematical Institute, Tohoku University, Sendai 980-8578, Japan} 
\email{daisuke.kazukawa.s6@dc.tohoku.ac.jp} 
\thanks{The author is supported by JSPS KAKENHI Grant Number 17J02121}
\date{September 26, 2019}

\begin{document}
\maketitle
\thispagestyle{empty}

\begin{abstract}
We investigate the relation between the concentration and the product of metric measure spaces. We have the natural question whether, for two concentrating sequences of metric measure spaces, the sequence of their product spaces also concentrates. A partial answer is mentioned in Gromov's book \cite{Grmv}. We obtain a complete answer for this question.
\end{abstract}

\tableofcontents

\section{Introduction}


In this paper, we study the relation between the concentration and the product of metric measure spaces. The concentration introduced by Gromov in \cite{Grmv} is defined as the convergence of metric measure spaces with respect to the observable distance, which is based on the theory of measure concentration phenomenon studied by L\'evy and V.~Milman. A typical example of the concentration is that the sequence $\{S^n(1)\}_{n \in \N}$ of $n$-dimensional unit spheres in $\R^{n+1}$ concentrates to a one-point metric measure space as $n \to \infty$, where the sphere $S^n(1)$ is endowed with the standard Riemannian metric and normalized volume measure, and $\N$ is the set of positive integers. It is one of the most important characteristics of concentration that such sequences of metric measure spaces whose dimensions are unbounded converge. A sequence of metric measure spaces concentrating to a one-point metric measure space is called a L\'evy family. Each of classical examples of the spaces that exhibit the measure concentration phenomenon corresponds to a L\'evy family.


In order to give a non-trivial example of sequences of metric measure spaces concentrating to a limit space which consists of at least two points, Gromov proved in \cite{Grmv}*{3.$\frac{1}{2}$.46.} that for a fixed metric measure space $X$ and a L\'evy family $\{Y_n\}_{n \in \N}$, the $l_p$-product space $X \times_p Y_n$, $p \in [1, + \infty]$, of $X$ and $Y_n$ concentrates to $X$ as $n \to \infty$. Given two metric measure spaces $X$, $Y$ and an extended real number $p \in [1, + \infty]$, the $l_p$-product space $X \times_p Y$ is defined as the Cartesian product set $X \times Y$ equipped with the $l_p$-metric $d_{l_p}$ and the product measure $m_X \otimes m_Y$. Gromov's argument implies the following conclusion: for a concentrating sequence $\{X_n\}_{n \in \N}$ of metric measure spaces and a L\'evy family $\{Y_n\}_{n \in \N}$, the sequence $\{X_n \times_p Y_n\}_{n \in \N}$ of their $l_p$-product spaces concentrates to the limit space of $\{X_n\}_{n \in \N}$. This fact gives us the following question: {\it for two sequences $\{X_n\}_{n \in \N}$ and $\{Y_n\}_{n \in \N}$ of metric measure spaces concentrating to metric measure spaces $X$ and $Y$ {\rm (}both are not necessarily a one-point space{\rm )} respectively, does the sequence $\{X_n \times_p Y_n\}_{n \in \N}$ of their $l_p$-product spaces concentrate to the $l_p$-product space $X \times_p Y$ of the limit spaces?} The main result in this paper gives an affirmative answer to this question.


In \cite{Grmv}, Gromov introduced not only the observable distance but also the box distance that induces a finer convergence than the concentration. For the box-convergence of metric measure spaces induced by the box distance, the above question is solved easily. That is, for two sequences $\{X_n\}_{n \in \N}$ and $\{Y_n\}_{n \in \N}$ of metric measure spaces box-converging to metric measure spaces $X$ and $Y$ respectively, the sequence $\{X_n \times_p Y_n\}_{n \in \N}$ of their $l_p$-product spaces box-converges to the $l_p$-product space $X \times_p Y$ of the limit spaces . Such a relation between the box-convergence and the product structure is summarized in Section 4.1 in this paper. Our main question is understood to be one of the questions whether the concentration also has a property that the box-convergence has.


Before describing the main theorem, we explain a product, which is a generalized notion of the $l_p$-product, of two metric measure spaces. Let $F \colon [0, +\infty)^2 \to [0, +\infty)$ be a continuous function satisfying the following condition: for any two metric spaces $(X, d_X)$ and $(Y, d_Y)$, the function
\begin{equation*}
d_F((x,y), (x',y')) := F(d_X(x, x'), d_Y(y, y'))
\end{equation*}
is a metric on $X \times Y$. Such a function $F$ is called a metric preserving function. Given two metric measure spaces $X$, $Y$ and a continuous metric preserving function $F \colon [0, +\infty)^2 \to [0, +\infty)$, the triple $(X \times Y, d_F, m_X \otimes m_Y)$ is a metric measure space. In this paper, this space is called the product space of $X$ and $Y$ generated by $F$ and is denoted by $X \times_F Y$. For any extended real number $p \in [1, +\infty]$, we define
\begin{equation*}
F_p(s, t) := \left\{\begin{array}{ll} (s^p + t^p)^{\frac{1}{p}} & \text{ if } p < +\infty, \\ \max{\{s, t\}} & \text{ if } p = +\infty.  \end{array}\right.
\end{equation*}
The function $F_p$ is a metric preserving function. The distance function $d_{F_p}$ accords with the $l_p$-metric $d_{l_p}$, and then the product space $X \times_{F_p} Y$ accords with the $l_p$-product space $X \times_p Y$. Thus, the product generated by the metric preserving functions is a generalization of the $l_p$-product. Other than the function $F_p$, for example, the two functions
\begin{equation*}
\log(e^s + e^t - 1), \quad \frac{1}{2}(s + t) + \frac{1}{2} \max{\{s, t\}}
\end{equation*}
are two of the easiest examples of metric preserving functions. However, general metric preserving functions are more complicated. We say that the function $F$ is an isotone if $F(s, t) \leq F(s', t')$ for all $s \leq s'$ and $t \leq t'$.  In general, such functions are not necessarily isotones. We describe some properties of such generalized product spaces in Section 3.1 of this paper, and show a lot more examples of metric preserving functions in Section 3.2.


The following main theorem gives an answer to the more general question than that stated before for the $l_p$-product.


\begin{thm}\label{main}
Let $F_n, F \colon [0, +\infty)^2 \to [0, +\infty)$, $n = 1, 2, \ldots$, be continuous metric preserving functions. Assume that $F_n$ converges pointwise to $F$ as $n \to \infty$. Then the following {\rm (1)} and {\rm (2)} are equivalent to each other.
\begin{enumerate}
\item \label{main1} For any two sequences $\{X_n\}_{n \in \N}$ and $\{Y_n\}_{n \in \N}$ of metric measure spaces concentrating to metric measure spaces $X$ and $Y$ respectively, the sequence $\{X_n \times_{F_n} Y_n \}_{n \in \N}$ of their product spaces concentrates to the product space $X \times_F Y$ as $n \to \infty$.
\item \label{main2} For any $s, t \in [0, +\infty)$, 
\begin{equation*}
\lim_{n\to\infty} (F_n(s, t) - \inf_{s \leq s'; \, t \leq t'} F_n(s', t')) = 0.
\end{equation*}
\end{enumerate}
\end{thm}

\begin{rem}\label{compcond}
We set
\begin{equation*}
I_n(s, t) : = F_n(s, t) - \inf_{s \leq s'; \, t \leq t'} F_n(s', t')
\end{equation*}
and consider the following five conditions.
\begin{enumerate}
\item The functions $F_n$ are isotones (i.e., $I_n \equiv 0$) for all $n \in \N$.
\item $\lim_{n\to\infty} \sup_{s, t \geq 0} I_n(s, t) = 0$.
\item $\lim_{n\to\infty} \sup_{0\leq s, t \leq D} I_n(s, t) = 0$ for any $D > 0$.
\item $\lim_{n\to\infty} I_n(s, t) = 0$ for any $s, t \in [0, +\infty)$.
\item The function $F$ is an isotone.
\end{enumerate}
Under the setting of Theorem \ref{main}, it is easy to see that $(1) \Rightarrow (2) \Rightarrow (3) \Rightarrow (4) \Rightarrow (5)$. On the other hand, we see that only $(4) \Rightarrow (3)$ holds among the converse implications (see Lemmas \ref{ptcpt}, \ref{equiv_main} and Example \ref{cex_compcond}). The condition (3) is also equivalent to the concentration of product spaces, that is, the condition (1) of Theorem \ref{main}. 
\end{rem}

Since the function $F_p$ is an isotone for all $p \in [1, +\infty]$, we obtain the following corollary. 

\begin{cor}\label{main_cor}
Let $\{X_n\}_{n \in \N}$ and $\{Y_n\}_{n \in \N}$ be two sequences of metric measure spaces concentrating to metric measure spaces $X$ and $Y$ respectively. Assume that $p_n \in [1, +\infty]$ converges to $p \in [1, +\infty]$ as $n \to \infty$. Then the sequence $\{X_n \times_{p_n} Y_n \}_{n \in \N}$ of their $l_{p_n}$-product spaces concentrates to the $l_{p}$-product space $X \times_p Y$ as $n \to \infty$.
\end{cor}


We prove Theorem \ref{main} in Section 4.2 and Section 5. In addition, from our main result, we construct new examples of sequences of metric measure spaces concentrating to a non-trivial limit space. Given a concentrating sequence $\{X_n\}_{n \in \N}$ of metric measure spaces, it had been not known even whether the sequence $\{X_n \times_p X_n\}_{n \in \N}$ of the $l_p$-product spaces concentrates. We describe in Section 4.3 a new specific example of the concentrating sequences that are obtained by applying our main result.

As another topic, the notion of metric preserving functions gives another question. Let $F \colon [0, +\infty) \to [0, +\infty)$ be a function such that for any metric space $(X, d_X)$, the function $F \circ d_X$ is a metric on $X$. This $F$ is also called a metric preserving function. Roughly speaking, such a function is a 1-dimensional version of a metric preserving function defined before. We consider the following question: {\it does a metric-transformed sequence of a concentrating sequence by some metric preserving functions concentrate too?} We obtain an answer, which is related with Theorem \ref{main}.

\begin{thm}\label{main_1dim}
Let $F_n, F \colon [0, +\infty) \to [0, +\infty)$, $n = 1, 2, \ldots$, be continuous metric preserving functions. Assume that $F_n$ converges pointwise to $F$ as $n \to \infty$. Then the following {\rm (1)} and {\rm (2)} are equivalent to each other.
\begin{enumerate}
\item \label{main_1dim1} For any sequence $\{X_n\}_{n \in \N}$ of metric measure spaces concentrating to a metric measure space $X$, the sequence $\{(X_n, F_n \circ d_{X_n}, m_X)\}_{n \in \N}$ concentrates to  $(X, F \circ d_{X}, m_X)$ as $n \to \infty$.
\item \label{main_1dim2} For any $s \in [0, +\infty)$, 
\begin{equation*}
\lim_{n \to \infty} (F_n(s) - \inf_{s \leq s'} F_n(s')) = 0.
\end{equation*}
\end{enumerate}
\end{thm}

The implication from (\ref{main_1dim2}) to (\ref{main_1dim1}) of Theorem \ref{main_1dim} is a corollary of Theorem \ref{main}. On the other hand, in Section 5, the proof of the converse implication of Theorem \ref{main_1dim} gives an essential idea to the proof of Theorem \ref{main}. We prove this theorem in Section 4.2 and Section 5 together with Theorem \ref{main}.

As a matter of fact, we are able to generalize Theorem \ref{main} to a statement for product spaces of $N$ metric measure spaces for any finite number $N$. This generalization is shown in Section 6 in this paper.

\begin{ack}
The author would like to thank Professor Takashi Shioya, Hiroki Nakajima, and Shinichiro Kobayashi for their comments and encouragement. He is also grateful to Professor Takumi Yokota and Ryunosuke Ozawa for their information and comments.
\end{ack}

\section{Preliminaries}

In this section, we describe the definitions and some properties of metric measure space, the box distance and the observable distance. We use most of these notions along \cite{MMG}. As for more details, we refer to \cite{MMG} and \cite{Grmv}*{Chapter 3.$\frac{1}{2}$}.

\subsection{Metric measure spaces}
Let $(X, d_X)$ be a complete separable metric space and $m_X$ a Borel probability measure on $X$. We call the triple $(X, d_X, m_X)$ a {\it metric measure space}, or an {\it mm-space} for short. We sometimes say that $X$ is an mm-space, in which case the metric and the measure of $X$ are respectively indicated by $d_X$ and $m_X$.

\begin{dfn}[mm-Isomorphism]
Two mm-spaces $X$ and $Y$ are said to be {\it mm-isomorphic} to each other if there exists an isometry $f \colon \supp{m_X} \to \supp{m_Y}$ such that $f_* m_X = m_Y$, where $f_* m_X$ is the push-forward measure of $m_X$ by $f$. Such an isometry $f$ is called an {\it mm-isomorphism}. Denote by $\mathcal{X}$ the set of mm-isomorphism classes of mm-spaces.
\end{dfn}

Note that an mm-space $X$ is mm-isomorphic to $(\supp{m_X}, d_X , m_X)$. We assume that an mm-space $X$ satisfies $X = \supp{m_X} $ unless otherwise stated.

\subsection{Observable diameter}

For a metric space $(X, d_X)$, we denote by $\Lip_1(X)$ the set of 1-Lipschitz functions on $X$.

The observable diameter is one of the most fundamental invariants of an mm-space.

\begin{dfn}[Partial and observable diameter]
Let $X$ be an mm-space. For a real number $\alpha \leq 1$, we define the {\it partial diameter} $\PD(X; \alpha)$ of $X$ to be the infimum of $\diam{A}$, where $A \subset X$ runs over all Borel subsets with $m_X(A) \geq \alpha$ and $\diam{A}$ is the diameter of $A$. For a real number $\kappa > 0$, we define the {\it observable diameter} of $X$ to be
\begin{equation}
\OD(X; -\kappa) := \sup_{f \in \Lip_1(X)} \PD((\R, |\cdot|, f_* m_X); 1 - \kappa).
\end{equation}
\end{dfn}

The observable diameter is an invariant under mm-isomorphism. Note that $\OD(X; -\kappa)$ is nonincreasing in $\kappa > 0$. 

\begin{dfn}[L\'evy family]
A sequence of mm-spaces $\{X_n\}_{n \in \N}$ is called a {\it L\'evy family} if
\begin{equation}
\lim_{n \to \infty} \OD(X_n; -\kappa) = 0
\end{equation}
for any $\kappa > 0$.
\end{dfn}

\subsection{Box distance and observable distance}
For a subset $A$ of a metric space $(X, d_X)$ and for a real number $r > 0$, we set
\begin{align}
U_r(A) & := \{x \in X \mid d_X(x, A) < r\}, 
\end{align}
where $d_X(x, A) := \inf_{a \in A} d_X(x, a)$.

\begin{dfn}[Prokhorov distance]
Let $\lambda > 0$ be a real number. The {\it $\lambda$-Prokhorov distance} $\prok^{(\lambda)}(\mu, \nu)$ between two Borel probability measures $\mu$ and $\nu$ on a metric space $X$ is defined to be the infimum of $\varepsilon > 0$ satisfying
\begin{equation}\label{prok}
\mu(U_\varepsilon(A)) \geq \nu(A) -\lambda \varepsilon
\end{equation}
for any Borel subset $A \subset X$. In particular, the $1$-Prokhorov distance $\prok^{(1)}$ is called the {\it Prokhorov distance} and we denote it by $\prok$.
\end{dfn}
The Prokhorov metric $\prok$ is a metrization of the weak convergence of Borel probability measures on $X$ provided that $X$ is a separable metric space.

\begin{dfn}[Ky Fan metric]
Let $(X, \mu)$ be a measure space and $(Y, d_Y)$ a metric space. For two $\mu$-measurable maps $f,g \colon X \to Y$, we define $\kf^\mu (f, g)$ to be the infimum of $\varepsilon \geq 0$ satisfying
\begin{equation}
\mu(\{x \in X \mid d_Y(f(x),g(x)) > \varepsilon \}) \leq \varepsilon.
\end{equation}
The two variable function $\kf^\mu$ is a metric on the set of $\mu$-measurable maps from $X$ to $Y$ by identifying two maps if they are equal to each other $\mu$-almost everywhere. We call $\kf^\mu$ the {\it Ky Fan metric}.
\end{dfn}

\begin{lem}[\cite{MMG}*{Lemma 1.26}]
Let $X$ be a topological space with a Borel probability measure $\mu$ and $Y$ a metric space. For any two Borel measurable maps $f, g \colon X \to Y$, we have
\begin{equation}
\prok(f_*\mu, g_* \mu) \leq \kf^\mu (f, g).
\end{equation}
\end{lem}

\begin{dfn}[Parameter]
Let $I := [0,1)$ and let $X$ be an mm-space. A map $\varphi \colon I \to X$ is called a {\it parameter} of $X$ if $\varphi$ is a Borel measurable map such that 
\begin{equation*}
\varphi_\ast \mathcal{L}^1 = m_X,
\end{equation*}
where $\mathcal{L}^1$ is the one-dimensional Lebesgue measure on $I$. 
\end{dfn}

\begin{lem}[\cite{MMG}*{Lemma 4.2}]\label{prm_exist}
Any mm-space has a parameter.
\end{lem}

\begin{dfn}[Box distance]
We define the {\it box distance} $\square(X, Y)$ between two mm-spaces $X$ and $Y$ to be the infimum of $\varepsilon \geq 0$ satisfying that there exist parameters $\varphi \colon I \to X$, $\psi \colon I \to Y$, and a Borel subset $I_0 \subset I$ with $\mathcal{L}^1(I_0) \geq 1 - \varepsilon$ such that
\begin{equation}
|d_X(\varphi(s), \varphi(t)) - d_Y(\psi(s), \psi(t))| \leq \varepsilon
\end{equation}
for any $s,t \in I_0$.
\end{dfn}

\begin{thm}[\cite{MMG}*{Theorem 4.10}]
The box distance function $\square$ is a complete separable metric on $\mathcal{X}$.
\end{thm}

\begin{lem}[\cite{MMG}*{Proposition 4.12}]\label{mmg4.12}
Let $X$ be a complete separable metric space. For any two Borel probability measures $\mu$ and $\nu$ on $X$, we have
\begin{equation}
\square((X, \mu), (X, \nu)) \leq 2\prok(\mu, \nu).
\end{equation}
\end{lem}

The following notion gives one of the conditions that are equivalent to the box convergence.

\begin{dfn}[$\varepsilon$-mm-Isomorphism]
Let $X$ and $Y$ be two mm-spaces and $f \colon X \to Y$ a Borel measurable map. Let $\varepsilon \geq 0$ be a real number. We say that $f$ is an {\it $\varepsilon$-mm-isomorphism} if there exists a Borel subset $X_0 \subset X$ such that 
\begin{enumerate}
\item $m_X(X_0) \geq 1 - \varepsilon$,
\item $|d_X(x, y) - d_Y(f(x),f(y))| \leq \varepsilon$ for any $x, y \in X_0$,
\item $\prok(f_* m_X, m_Y) \leq \varepsilon$.
\end{enumerate}
We call $X_0$ a {\it nonexceptional domain} of $f$.
\end{dfn}
It is easy to see that, for a 0-mm-isomorphism $f \colon X \to Y$, there is an mm-isomorphism $\hat{f} \colon X \to Y$ that is equal to $f$ $m_X$-a.e.~on $X$.

\begin{lem}[\cite{MMG}*{Lemma 4.22}]
\begin{enumerate}
\item[]
\item If there exists an $\varepsilon$-mm-isomorphism $f \colon X \to Y$, then $\square(X, Y) \leq 3\varepsilon$.
\item If $\square(X, Y) < \varepsilon$, then there exists a $3\varepsilon$-mm-isomorphism $f\colon X \to Y$.
\end{enumerate}
\end{lem}

For any topological space $X$, any metric space $Y$, and any Borel measurable map $p \colon X \to Y$, we set
\begin{equation}\label{pull_lip1}
p^* \Lip_1(Y) := \{ f \circ p \mid f \in \Lip_1(Y) \}.
\end{equation}
Note that, for any mm-space $X$ and any parameter $\varphi \colon I \to X$ of $X$, the set $\varphi^* \Lip_1(X)$ consists of Borel measurable functions on $I$.

\begin{dfn}[Observable distance] 
We define the {\it observable distance} $\conc(X, Y)$ between two mm-spaces $X$ and $Y$ by
\begin{equation*}
\conc(X, Y) := \inf_{\varphi, \psi} \haus(\varphi^* \Lip_1(X), \psi^* \Lip_1(Y)),
\end{equation*}
where $\varphi \colon I \to X$ and $\psi \colon I \to Y$ run over all parameters of $X$ and $Y$ respectively, and $\haus$ is the Hausdorff distance with respect to the metric $\kf^{\mathcal{L}^1}$. We say that a sequence of mm-spaces $\{X_n\}_{n \in \N}$ {\it concentrates} to an mm-space $X$ if $X_n$ $\conc$-converges to $X$ as $n \to \infty$.
\end{dfn}

\begin{prop}[\cite{MMG}*{Corollary 5.8}]
Let $\{X_n\}_{n \in \N}$ be a sequence of mm-spaces. Then, $\{X_n\}_{n \in \N}$ is a L\'evy family if and only if $\{X_n\}_{n \in \N}$ concentrates to a one-point mm-space as $n \to \infty$.
\end{prop}

\begin{ex}[\cite{MMG}*{Section 2.3}]\label{sphere_riem}
Let $S^n(r_n)$, $n = 1, 2, \ldots$, be the sphere of radius $r_n > 0$ in $\R^{n+1}$. Assume that $S^n(r_n)$ endowed with the standard Riemannian metric. Let $\sigma^n$ be the Riemannian volume measure on $S^n(r_n)$ normalized as $\sigma^n(S^n(r_n)) = 1$. Then we have
\begin{equation*}
\OD((S^n(r_n), d_{S^n(r_n)}, \sigma^n); -\kappa) = O(r_n \, n^{-1/2})
\end{equation*}
for any $\kappa > 0$ as $n \to \infty$. That is, the sequence $\{S^n(r_n)\}_{n \in \N}$ is a L\'evy family if $r_n = o(\sqrt{n})$.
\end{ex}

\begin{prop}[\cite{MMG}*{Proposition 5.5}]
For any two mm-spaces $X$ and $Y$, we have
\begin{equation}
\conc(X, Y) \leq \square(X, Y).
\end{equation}
\end{prop}

\begin{dfn}[Enforce $\varepsilon$-concentration]
A Borel measurable map $p \colon X \to Y$ is said to {\it enforce $\varepsilon$-concentration of $X$ to $Y$} if 
\begin{equation}
\haus(\Lip_1(X), p^* \Lip_1(Y)) \leq \varepsilon,
\end{equation}
where $\haus$ is the Hausdorff distance with respect to the metric $\kf^{m_X}$. 
\end{dfn}

\begin{thm}[\cite{MMG}*{Corollary 5.36}]\label{equiconc}
Let $X_n$ and $X$ be mm-spaces, where $n = 1, 2, \ldots$. Then the following {\rm (1)} and {\rm (2)} are equivalent to each other.
\begin{enumerate}
\item $\{X_n\}_{n \in \N}$ concentrates to $X$ as $n \to \infty$.
\item There exists a sequence of Borel measurable maps $p_n \colon X_n \to X$, $n = 1, 2, \ldots$, that enforce $\varepsilon_n$-concentration of $X_n$ to $X$ and $\prok((p_n)_* m_{X_n}, m_X) \leq \varepsilon_n$ for all $n$ and for some sequence $\varepsilon_n \to 0$.
\end{enumerate}
\end{thm}

\subsection{Strassen's theorem and L\'evy mean}
\begin{dfn}[Transport plan]
Let $\mu$ and $\nu$ be two finite Borel measures on $X$. A Borel measure $\pi$ on $X \times X$ is called a {\it transport plan} (or {\it coupling}) between $\mu$ and $\nu$ if
\begin{equation}
\pi(A \times X) = \mu(A) \text{ \ and \ } \pi(X \times A) = \nu(A) 
\end{equation}
for any Borel subset $A \subset X$.
\end{dfn}

\begin{dfn}[$\varepsilon$-subtransport plan]
Let $\mu$ and $\nu$ be two Borel probability measures on $X$. A Borel measure $\pi$ on $X \times X$ is called an {\it $\varepsilon$-subtransport plan} between $\mu$ and $\nu$ if there exist two Borel measures $\mu'$ and $\nu'$ on $X$ with $\mu' \leq \mu$ and $\nu' \leq \nu$ such that $\pi$ is a transport plan between $\mu'$ and $\nu'$, and $\pi$ satisfies
\begin{equation}
\supp{\pi} \subset \Delta_\varepsilon := \{ (x, x') \in X \times X \mid d_X(x, x') \leq \varepsilon \}.
\end{equation}
For an $\varepsilon$-subtransport plan $\pi$ between $\mu$ and $\nu$, the {\it deficiency} of $\pi$ is defined to be
\begin{equation}
\df{\pi} := 1 - \pi(X \times X).
\end{equation}
\end{dfn}

\begin{thm}[Strassen's theorem]\label{Strassen}
Assume that $X$ is a complete separable metric space. For any real number $\lambda > 0$ and for any two Borel probability measures $\mu$ and $\nu$ on $X$, we have
\begin{equation}
\lprok{\lambda}(\mu, \nu) = \inf{\left\{ \varepsilon > 0 \midd \begin{array}{l} \text{There exists an } \varepsilon\text{-subtransport plan } \pi \\ \text{between } \mu \text{ and } \nu \text{ with } \df{\pi} \leq \lambda\varepsilon \end{array} \right\}}.
\end{equation}
\end{thm}

\begin{dfn}[Median and L\'evy mean]
Let $X$ be a measure space with probability measure $\mu$ and $f \colon X \to \R$ a measurable function. A real number $m \in \R$ is called a {\it median} of $f$ if it satisfies
\begin{equation*}
\mu(\{ x \in X \mid  f(x) \geq m\}) \geq \frac{1}{2} \text{ \ and \ } \mu(\{ x \in X \mid  f(x) \leq m\}) \geq \frac{1}{2}.
\end{equation*}
It is easy to see that the set of medians of $f$ is a closed and bounded interval. The {\it L\'evy mean} $\lm(f; \mu)$ of $f$ with respect to $\mu$ is defined to be
\begin{equation}
\lm(f; \mu) := \frac{\underline{m} +\overline{m}}{2},
\end{equation}
where $\underline{m}$ is the minimum of medians of $f$, and $\overline{m}$ the maximum of medians of $f$.
\end{dfn}

\begin{prop}[\cite{MMG}*{Section 2.3}]
Let $\{X_n\}_{n \in \N}$ be a sequence of mm-spaces. Then, $\{X_n\}_{n \in \N}$ is a L\'evy family if and only if for any $f_n \in \Lip_1(X_n)$, 
\begin{equation}
\lim_{n \to \infty} \kf^{m_{X_n}}(f_n, \lm(f_n; m_{X_n})) = 0.
\end{equation}
\end{prop}

\begin{dfn}[L\'evy radius]
Let $X$ be an mm-space and $\kappa > 0$ a real number. The {\it L\'evy radius} $\LR(X; -\kappa)$ of $X$ is defined to be the infimum of $\varepsilon > 0$ satisfying
\begin{equation}
m_X(\{ x \in X \mid | f(x) - \lm(f; m_X) | > \varepsilon\}) \leq \kappa
\end{equation}
for any $f \in \Lip_1(X)$.
\end{dfn}

\begin{lem}[\cite{MMG}*{Lemma 7.31}]\label{mmg7.31}
Let $X$ be an mm-space. For any $\kappa$ with $0 < \kappa < 1/2$, we have
\begin{equation}
\LR(X; -\kappa) \leq \OD(X; - \kappa).
\end{equation}
\end{lem}

\begin{lem}[\cite{MMG}*{Lemma 9.6}]\label{lm_lem}
Let $\mu$ and $\nu$ be two Borel probability measures on a metric space $X$. Assume that there exists an $\varepsilon$-subtransport plan $\pi$ between $\mu$ and $\nu$ with $\df{\pi} < 1-2\kappa$ for two real numbers $\varepsilon$ and $\kappa$ with $\varepsilon > 0$ and $0 < \kappa < 1/2$. Then, for any 1-Lipschitz function $f \colon X \to \R$, we have
\begin{equation}
\begin{split}
& |\lm(f ; \mu) - \lm(f ; \nu)| \\
& \leq \varepsilon + \OD((X, \mu) ; -\kappa) + \OD((X, \nu) ; -\kappa).
\end{split}
\end{equation}
\end{lem}

\section{Product space of metric measure spaces}

\subsection{Metric preserving functions and product spaces}

\begin{dfn}[Metric preserving function]
Let $N \in \N$. A function $F \colon [0, +\infty)^N \to [0, +\infty)$ is called a {\it metric preserving function} if for any $N$ metric spaces $(X_1, d_{X_1}), \ldots, (X_N, d_{X_N})$, the function
\begin{equation}\label{defdF}
d_F((x_i)_{i=1}^N, (x'_i)_{i=1}^N) := F(d_{X_1}(x_1, x'_1), \ldots, d_{X_N}(x_N, x'_N))
\end{equation}
is a metric on $X_1 \times \cdots \times X_N$.
\end{dfn}

Note that for a metric preserving function $F$,
\begin{equation*}
F^{-1}(0) = \{(0, \ldots, 0)\}
\end{equation*}
holds necessarily. Let $a$, $b$, and $c$ be three nonnegative real numbers. We call the triplet $(a, b, c)$ a {\it triangle triplet} if $a \leq b + c$, $b \leq a + c$, and $c \leq a + b$ are all satisfied.

\begin{thm}[\cite{BD}*{Theorem 2.6}]\label{TT}
A function $F \colon [0, +\infty)^N \to [0, +\infty)$ with $F^{-1}(0) = \{(0,  \ldots, 0)\}$ is a metric preserving function if and only if for any $N$ triangle triplets $(a_i, b_i, c_i)$, $i=1, \ldots, N$, 
the triplet 
\begin{equation*}
(F(a_1, \ldots, a_N), F(b_1, \ldots, b_N), F(c_1, \ldots, c_N))
\end{equation*}
is a triangle triplet.
\end{thm}

\begin{cor}\label{cor_mpf}
Let $F \colon [0, +\infty)^N \to [0, +\infty)$ be a metric preserving function. Then, for any $(s_i)_{i=1}^N, (s'_i)_{i=1}^N \in [0, +\infty)^N$,
\begin{enumerate}
\item $|F(s_1, \ldots, s_N) - F(s'_1, \ldots, s'_N)| \leq F(|s_1-s'_1|, \ldots |s_N-s'_N|)$,
\item $F(s_1, \ldots, s_N) \leq 2 F(s'_1, \ldots, s'_N)$ if $s_i \leq 2 s'_i$ for every $i$.
\end{enumerate}
\end{cor}

\begin{proof}
For any $s, s' \in [0, +\infty)$, the triplet $(s, s', |s - s'|)$ is a triangle triplet, and the triplet $(s, s', s')$ is a triangle triplet if $s \leq 2s'$. Applying Theorem \ref{TT} to them, we obtain this corollary.
\end{proof}

\begin{dfn}
Let $N \in \N$. A function $F \colon [0, +\infty)^N \to [0, +\infty)$ is said to be {\it subadditive} if for any $(s_i)_{i=1}^N, (s'_i)_{i=1}^N \in [0, +\infty)^N$, 
\begin{equation*}
F(s_1 + s'_1, \ldots, s_N + s'_N) \leq F(s_1, \ldots, s_N)  + F(s'_1, \ldots, s'_N).
\end{equation*}
$F$ is called an {\it isotone} if 
\begin{equation*}
F(s_1, \ldots, s_N) \leq F(s'_1, \ldots, s'_N)
\end{equation*}
for any $(s_i)_{i=1}^N, (s'_i)_{i=1}^N \in [0, +\infty)^N$ such that $s_i \leq s'_i$ for each $i$,
\end{dfn}

In the case of $N = 1$, an isotone means a nondecreasing function.

\begin{lem}[\cite{Kelley}*{Exercise 4.C}, \cite{Tibor}*{Satz 1}]\label{Kelly}
Let $N \in \N$ and let $F \colon [0, +\infty)^N \to [0, +\infty)$ be a function. Then the following {\rm (1)} and {\rm (2)} hold.
\begin{enumerate}
\item If $F$ is a metric preserving function, then $F$ is subadditive.
\item If $F$ is a subadditive isotone fulfilling $F^{-1}(0) = \{(0, \ldots, 0)\}$, then $F$ is a metric preserving function.
\end{enumerate}
\end{lem}

\begin{rem}
There are some examples of metric preserving functions that are not isotones even if $N = 1$. We show such examples in the next subsection.
\end{rem}

In this paper, we usually consider the case of $N=2$.

\begin{ex}\label{exF}
We cite some easy (however important) examples of metric preserving functions.
\begin{enumerate}
\item\label{exF1} 
For any extended real number $p \in [1, +\infty]$, we define
\begin{equation*}
F_p(s, t) := \left\{\begin{array}{ll} (s^p + t^p)^{\frac{1}{p}} & \text{ if } p < +\infty, \\ \max{\{s, t\}} & \text{ if } p = +\infty.  \end{array}\right.
\end{equation*}
\item\label{exF2} Let $F_{\exp}(s, t) := \log(e^s + e^t -1)$.
\item\label{exF3} 
For any real number $\alpha \in (0, 1)$, we define
\begin{equation*}
F_\alpha(s, t) := s^\alpha + t^\alpha.
\end{equation*}
\item\label{exF4} 
For two real numbers $p, q \in [1, +\infty)$ with $p < q$, we define
\begin{equation*}
F_{p, q}(s, t) := (s^p + t^p)^{\frac{1}{q}}.
\end{equation*}
\end{enumerate}
\end{ex}


\begin{thm}[\cite{BD}*{Section 3}]
Let $F \colon [0, +\infty)^N \to [0, +\infty)$ be a metric preserving function. Then the following conditions are equivalent to each other.
\begin{enumerate}
\item $F$ is continuous.
\item $F$ is continuous at $(0, \ldots, 0)$.
\item $F$ is uniformly continuous.
\item For each $i = 1, \ldots, N$, the map $F \circ \iota_i \colon [0, +\infty) \to [0, +\infty)$ is continuous, where $\iota_i \colon [0, +\infty) \to [0, +\infty)^N$ is the natural $i$-th inclusion map {\rm (}i.e., $\pr_i \circ \iota_i = \id$ and $\pr_j \circ \iota_i = 0$ for $j \neq i${\rm )}.
\item For any $N$ metric spaces $(X_1, d_{X_1}), \ldots, (X_N, d_{X_N})$, the metric $d_F$ induces the product topology on $X_1 \times \cdots \times X_N$.
\end{enumerate}
\end{thm}

\begin{rem}
In the case of $N = 1$, we understand the condition (5) in the above theorem to mean that $(X, F \circ d_X)$ has the same topology as $(X, d_X)$.
\end{rem}

\begin{prop}
Let $F \colon [0, +\infty)^N \to [0, +\infty)$ be a continuous metric preserving function. If $N$ metric spaces $X_1, \ldots, X_N$ are complete, then the metric space $(X_1 \times \cdots \times X_N, d_F)$ is also complete.
\end{prop}

\begin{proof}
We take any $d_F$-Cauchy sequence $\{(x_n^1, \ldots, x_n^N)\}_{n \in \N}$. We have 
\begin{equation*}
\lim_{m,n \to \infty} F(d_{X_1}(x_m^1, x_n^1), \ldots, d_{X_N}(x_m^N, x_n^N)) = 0.
\end{equation*}
Fix $i \in \{1, \ldots, N\}$. Since
\begin{align*}
& F(0, \ldots, 0, d_{X_i}(x_m^i, x_n^i), 0, \ldots , 0) \\
& \leq 2F(d_{X_1}(x_m^1, x_n^1), \ldots, d_{X_N}(x_m^N, x_n^N)) \to 0
\end{align*}
as $m, n \to \infty$, Corollary \ref{cor_mpf} (2) leads to
\begin{equation*}
\lim_{m,n \to \infty} d_{X_i}(x_m^i, x_n^i) = 0.
\end{equation*}
By the completeness of $X_i$, there exists $x^i \in X_i$ such that $x_n^i \to x^i$ in $X_i$ as $n \to \infty$. Thus
\begin{align*}
& \lim_{n \to \infty} d_F((x_n^1, \ldots, x_n^N), (x^1, \ldots, x^N)) \\
& = \lim_{n \to \infty} F(d_{X_1}(x_n^1, x^1), \ldots, d_{X_N}(x_n^N, x^N)) = 0,
\end{align*}
which implies that $d_F$ is complete.
\end{proof}



We denote by $\F^N$ the set of continuous metric preserving functions on $[0, +\infty)^N$.

\begin{dfn}[Product space]
Given two mm-spaces $X$, $Y$ and a function $F \in \F^2$, we define the {\it product space} $X \times_F Y$ of $X$ and $Y$ as the mm-space 
\begin{equation}
(X \times Y, d_F, m_X \otimes m_Y)
\end{equation}
which consists of the complete separable metric space $(X \times Y, d_F)$ and the product measure $m_X \otimes m_Y$ of $m_X$ and $m_Y$.
\end{dfn}

\begin{dfn}[$l_p$-Product space]
For two mm-spaces $X$, $Y$ and an extended real number $p \in [1, +\infty]$, we call the distance function $d_{F_p}$ on $X \times Y$, where $F_p$ is of Example \ref{exF} (\ref{exF1}), the {\it $l_p$-metric}, and the product space $X \times_{F_p} Y$ the {\it $l_p$-product space}. From now on, we denote $X \times_{F_p} Y$ by $X \times_p Y$.
\end{dfn}

\subsection{Examples of metric preserving function}


In this subsection, we present many examples of continuous metric preserving functions. At first, we are able to make such functions from the Mulholland inequality known as the generalized Minkowski inequality.



\begin{thm}[Mulholland inequality \cite{Mul}, \cite{Kuc}*{Theorem 8.8.1}]\label{Mul}
If both a homeomorphism $\varphi \colon [0, +\infty) \to [0, +\infty)$ {\rm(}i.e., continuous and increasing bijection with $\varphi(0) = 0${\rm )} and $\log \circ \varphi \circ \exp$ are convex, then
\begin{equation}\label{Mulineq}
\varphi^{-1}(\varphi(s + s') + \varphi(t + t')) \leq \varphi^{-1}(\varphi(s) + \varphi(t)) + \varphi^{-1}(\varphi(s') + \varphi(t'))
\end{equation}
for any $s, s', t, t' \in [0, +\infty)$.
\end{thm}

For a homeomorphism $\varphi \colon [0, +\infty) \to [0, +\infty)$, the function $F_\varphi$ on $[0, +\infty)^2$ defined by
\begin{equation}\label{genfct}
F_{\varphi}(s, t) := \varphi^{-1}(\varphi(s) + \varphi(t))
\end{equation}
for $s, t \in [0, +\infty)$ is an isotone with $F^{-1}(0) = \{(0, 0)\}$. Mulholland inequality says that if both $\varphi$ and $\log \circ \varphi \circ \exp$ are convex in addition, then $F_{\varphi}$ is subadditive, so that it is a metric preserving function.

\begin{lem}[\cite{Kuc}*{Lemma 8.8.1}]
We consider the following two classes $A$ and $B$ of functions.
\begin{enumerate}
\item[(A)] The class $A$ consists of homeomorphisms $\varphi \colon [0, +\infty) \to [0, +\infty)$ such that both $\varphi$ and $\log \circ \varphi \circ \exp$ are convex.
\item[(B)] The class $B$ consists of the functions $\gamma \colon \R \to \R$ that are continuous, increasing, and convex.
\end{enumerate}
For any function $\varphi$ in $A$, we define a function $\gamma$ belonging to $B$ as
\begin{equation}
\gamma := \log \circ \varphi \circ \exp - \id_{\R}
\end{equation}
Then this gives a one to one correspondence between $A$ and $B$.
\end{lem}

\begin{ex}
We cite below some specific examples of functions fulfilling the assumptions of Theorem \ref{Mul}.

\begin{enumerate}
\item 
For a real number $p \in [1, +\infty)$, let $\varphi(s) := s^p$. Then the function $F_\varphi$ of (\ref{genfct}) coincides with $F_p$ in Example \ref{exF} (\ref{exF1}).
\item Let $\varphi(s) := e^s-1$. Then the function $F_\varphi$ of (\ref{genfct}) coincides with $F_{\exp}$ in Example \ref{exF} (\ref{exF2}).
\item Let $\varphi(s) := \sinh{s}$. Then $F_\varphi(s, t) = \arcsinh{(\sinh{s} + \sinh{t})}$.
\item Let $\varphi(s) := s^2 + 2s$. Then $F_\varphi(s, t) = \sqrt{s^2+2s+t^2+2t+1}-1$.
\end{enumerate}
\end{ex}

\begin{rem}
In \cite{Mul}, Mulholland proved that if a homeomorphism $\varphi \colon [0, +\infty) \to [0, +\infty)$ satisfies (\ref{Mulineq}), then $\varphi$ is convex. In recent years, Petr\'ik \cite{Pet} gave new type examples of functions that fulfill \eqref{Mulineq}. One of his specific examples is
\begin{equation}
\varphi(s) := \left\{ \begin{array}{ll} \frac{5}{3}s & \text{ if } s \in [0, 1), \\ \frac{7}{3}s - \frac{2}{3} & \text{ if } s \in [1, 2), \\ s^2 & \text{ if } s \in [2, +\infty). \end{array} \right.
\end{equation}
This $\varphi$ satisfies (\ref{Mulineq}) but $\log \circ \varphi \circ \exp$ is not convex. Of course, for this $\varphi$, the function $F_\varphi$ of (\ref{genfct}) is a metric preserving function.
\end{rem}


\begin{prop}\label{comb}
Let $F_1$, $F_2$, and $F$ be metric preserving functions on $[0, +\infty)^2$ and let $f_1$, $f_2$, and $f$ metric preserving functions on $[0, +\infty)$. Then, the functions $G_1$, $G_2$, and $G_3$ defined as
\begin{align*}
G_1(s, t) & := F_1(s, t) + F_2(s, t), \\
G_2(s, t) & := f_1(s) + f_2(t),  \\
G_3(s, t) & := f(F(f_1(s), f_2(t)))
\end{align*}
for $s, t \in [0, +\infty)$ are metric preserving functions on $[0, +\infty)^2$. 
\end{prop}

\begin{proof}
This proposition follows from Theorem \ref{TT} directly.
\end{proof}

\begin{rem}
It is possible to make the functions $F_\alpha$ and $F_{p, q}$ in Example \ref{exF} applying Proposition \ref{comb} to the function $F_p$ in Example \ref{exF} (\ref{exF1}).
Applying Proposition \ref{comb} to $F_1(s, t) = s + t$ with $f_1(s) = f_2(s) = s^\alpha$, $\alpha \in (0, 1]$, and $g = \id$, we obtain $F_\alpha(s, t) = s^\alpha + t^\alpha$ in Example \ref{exF} (\ref{exF3}). Applying Proposition \ref{comb} to $F_p(s, t) = (s^p + t^p)^{1/p}$, $p \in [1, +\infty)$ with $f_1 = f_2 = \id$ and $g(s) = s^{p/q}$, $q \in [1, +\infty)$ such that $p < q$, we obtain $F_{p,q}(s, t) = (s^{p} + t^{p})^{1/q}$ in Example \ref{exF} (\ref{exF4}).
\end{rem}

We show some examples of metric preserving functions that are not isotones. We say that a function $F \colon [0, +\infty)^N \to [0, +\infty)$ is {\it tightly bounded} if there exists $T > 0$ such that $T \leq F(s_1, \ldots, s_N) \leq 2T$ for every $(s_i)_{i=1}^N \in [0, +\infty)^N \setminus \{(0, \ldots, 0)\}$.

\begin{prop}
If a function $F \colon [0, +\infty)^N \to [0, +\infty)$ is tightly bounded and $F(0, \ldots, 0) = 0$, then $F$ is a metric preserving function {\rm (}but is not continuous{\rm )}.
\end{prop}

\begin{proof}
We take any $N$ triangle triplets $\{(a_i, b_i, c_i)\}_{i=1}^N$. By symmetry, it is sufficient to prove
\begin{equation}\label{tbprop_goal}
F(a_1, \ldots, a_N) \leq F(b_1, \ldots, b_N) + F(c_1, \ldots, c_N).
\end{equation}
If $(b_1, \ldots, b_N) = (0, \ldots, 0)$, then $(a_1, \ldots, a_N) = (c_1, \ldots, c_N)$, which implies \eqref{tbprop_goal}. In the case of $(c_1, \ldots, c_N) = (0, \ldots, 0)$, we have \eqref{tbprop_goal} similarly. If $(b_1, \ldots, b_N), (c_1, \ldots, c_N) \neq (0, \ldots, 0)$, then
\begin{equation*}
F(a_1, \ldots, a_N) \leq 2T = T + T \leq F(b_1, \ldots, b_N) + F(c_1, \ldots, c_N)
\end{equation*}
for some $T > 0$ since $F$ is tightly bounded. This completes the proof.
\end{proof} 

\begin{prop}[\cite{DP}*{Theorem 1}]\label{DP}
Let $F$ and $G$ be two metric preserving functions on $[0, +\infty)$. Assume that there exists $r > 0$ such that $F(r) = G(r)$. Define  a function $H$ on $[0, + \infty)$ by
\begin{equation*}
H(s) := \left\{ \begin{array}{ll} F(s) & \text{ if } s \in [0, r), \\ G(s) & \text{ if } s \in [r, +\infty). \end{array} \right.
\end{equation*}
If $F$ is concave and $|G(s) - G(t)| \leq F(|s - t|)$ holds for any $s, t \in [r, +\infty)$ with $|s - t| \leq r$, then $H$ is a metric preserving function.
\end{prop}

\begin{ex}
Proposition \ref{DP} gives the following specific examples of continuous metric preserving functions that are not isotones.
\begin{align*}
H_1(s) & := \left\{ \begin{array}{ll} s & \text{ if } s \in [0, 2), \\ 4 - s & \text{ if } s \in [2, 3), \\ 1 & \text{ if } s \in [3, +\infty). \end{array} \right. \\
H_2(s) & := \left\{ \begin{array}{ll} s & \text{ if } s \in [0, 1), \\ \displaystyle \frac{1+s+\sin^2(s-1)}{2s} & \text{ if } s \in [1, +\infty). \end{array} \right.
\end{align*}
\end{ex}

Given a function $F \in \F^1$ that is not an isotone, the function $G$ on $[0, +\infty)^2$ defined by
\begin{equation*}
G(s, t) := F(s) + F(t)
\end{equation*}
belongs to $\F^2$ (by Proposition \ref{comb}) but is not an isotone. The following example is related with Remark \ref{compcond}.

\begin{ex}\label{cex_compcond}
We define three functions $F_n^1$, $F_n^2$, and $F_n^3$ by
\begin{align*}
F_n^1(s) & := \left\{ \begin{array}{ll} s & \text{ if } s \in [0, 2), \\  4 - s & \text{ if } s \in [2, 2+n^{-1}), \\ 2 - n^{-1} & \text{ if } s \in [2 + n^{-1}, +\infty). \end{array} \right. \\
F_n^2(s) & := \left\{ \begin{array}{ll} s & \text{ if } s \in [0, 2), \\ 2 & \text{ if } s \in [2, n+2), \\  s - n & \text{ if } s \in [n+2, n+3), \\ n + 6 - s & \text{ if } s \in [n+3, n+4), \\ 2 & \text{ if } s \in [n+4, +\infty). \end{array} \right. \\
F_n^3(s) & := \left\{ \begin{array}{ll} s & \text{ if } s \in [0, 2), \\ 2 & \text{ if } s \in [2, n+2), \\  n + 4 - s & \text{ if } s \in [n+2, n+3), \\ 1 & \text{ if } s \in [n+3, +\infty). \end{array} \right.
\end{align*}
By Proposition \ref{DP}, the functions $F_n^i$, $i = 1, 2, 3$, belong to $\F^1$ and converge to the function $\min\{s, 2\}$ as $n \to \infty$. Let
\begin{equation*}
G_n^i(s, t) := F_n^i(s) + F_n^i(t)
\end{equation*}
for every $i = 1, 2, 3$. Note that $G_n^i$ belongs to $\F^2$ and converges to $\min\{s, 2\} + \min\{t, 2\}$ for every $i$. The functions $G_n^1$, $G_n^2$, and $G_n^3$ are counterexamples of ``(2) $\Rightarrow$ (1)'', ``(3) $\Rightarrow$ (2)'', and ``(5) $\Rightarrow$ (4)'' of Remark \ref{compcond} respectively.
\end{ex}

\begin{ex}\label{notisotone}
The following function $F$ is a continuous metric preserving function such that $F$ is not an isotone but both $s \mapsto F(s, 0)$ and $t \mapsto F(0, t)$ are nondecreasing.
\begin{equation*}
F(s, t) := \left\{ \begin{array}{ll} \min\{s, 1\} + \min\{t, 1\} & \text{ if } s \in [0, 1) \text{ or } t \in [0, 1), \\ 2 - \min\{s - 1, t - 1, 1 \}& \text{ if } s, t \in [1, +\infty). \end{array} \right.
\end{equation*}
\end{ex}

\subsection{Estimates of observable diameter}
In this subsection, we study the relation between the observable diameter and the metric preserving function. Our goal is the estimate of the observable diameter of product spaces.

\begin{dfn}[Concentration function]
Let $X$ be an mm-space. We define the {\it concentration function} $\alpha_X$ of $X$ to be
\begin{equation}
\alpha_{X}(r) := \sup_A (1 - m_X(U_r(A)))
\end{equation}
for $r > 0$, where $A \subset X$ runs over all Borel subsets with $m_X(A) \geq 1/2$.
\end{dfn}

\begin{prop}[\cite{Led}*{Proposition 1.12}, \cite{MMG}*{Remark 2.28}]\label{conc_fct}
\
\begin{enumerate}
\item $\OD(X; -\kappa) \leq 2\inf\{r > 0 \, | \, \alpha_X(r) \leq \kappa/2\}$ for any $\kappa > 0$.
\item $\alpha_X(r) \leq \sup\left\{\kappa >0  \midd \OD(X; -\kappa) \geq r\right\}$ for any $r > 0$.
\end{enumerate}
\end{prop}

\begin{lem}\label{key_1dim}
Let $F \in \F^1$ and let $X$ be an mm-space. Denote the mm-space $(X, F \circ d_X, m_X)$ by $X^F$. Then, we have 
\begin{enumerate}
\item $\sup_{\varepsilon>0} \alpha_{X^F}(2F(s) + \varepsilon) \leq \alpha_X(s)$ for any $s > 0$,
\item $\OD(X^F ; -2\kappa) \leq 4 F(\OD(X; -\kappa))$ for any $\kappa > 0$.
\end{enumerate}
\end{lem}

\begin{proof}
We first prove (1). Let $\varepsilon$ and $s$ be arbitrary positive real numbers. We take any $r > 0$ such that $r < \alpha_{X^F}(2F(s) + \varepsilon)$. There exists a Borel subset $A \subset X$ such that $m_X(A) \geq 1/2$ and 
\begin{equation*}
1- m_X(U_{2F(s)+\varepsilon}^{F}(A)) > r,
\end{equation*}
where $U_r^F(A)$ is the $r$-open neighborhood of $A$ with respect to $F \circ d_X$. By Corollary \ref{cor_mpf} (2), we have $U_s(A) \subset U_{2F(s)+\varepsilon}^{F}(A)$. In fact, for $x \in X$, if there exists $x' \in A$ such that $d_X(x, x') < s$, then $F(d_X(x, x')) \leq 2F(s)$ by Corollary \ref{cor_mpf} (2). We obtain 
\begin{equation*}
\inf_{x' \in A} F(d_X(x, x')) < 2F(s) + \varepsilon,
\end{equation*}
which implies $U_s(A) \subset U_{2F(s)+\varepsilon}^{F}(A)$. By this, we have 
\begin{equation*}
r < 1 - m_X(U_{2F(s)+\varepsilon}^{F}(A)) \leq 1 - m_X(U_s(A)) \leq \alpha_X(s).
\end{equation*}
Since $r$ is arbitrary, we obtain $\alpha_{X^F}(2F(s) + \varepsilon) \leq \alpha_X(s)$.

We next prove (2). Let $\varepsilon$ and $\kappa$ be arbitrary positive real numbers. Setting
\begin{equation*}
s_\varepsilon := \OD(X; -\kappa) + \varepsilon,
\end{equation*}
by (1) of this proposition and Proposition \ref{conc_fct} (2), we have 
\begin{equation*}
\alpha_{X^F}(2F(s_\varepsilon) + \varepsilon) \leq \alpha_X(s_\varepsilon) \leq \kappa.
\end{equation*}
Combining this and Proposition \ref{conc_fct} (1) implies
\begin{equation*}
\OD(X^F; -2\kappa) \leq 4F(s_\varepsilon) + 2\varepsilon.
\end{equation*}
As $\varepsilon \to 0$, we obtain $\OD(X^F ; -2\kappa) \leq 4 F(\OD(X; -\kappa))$. The proof is completed.
\end{proof}

\begin{lem}[\cite{MMG}*{in the proof of Proposition 7.32}]\label{lm_lip}
Let $X$ and $Y$ be two mm-spaces and let $p \in [1, +\infty]$. Given $f \in \Lip_1(X \times_p Y)$, we define functions $g \colon X \to \R$ and $h \colon Y \to \R$ by
\begin{equation}\label{lm_lip_eq}
g(x) := \lm(f(x, \cdot) ; m_Y) \text{ and } h(y) := \lm(f(\cdot, y) ; m_X)
\end{equation}
for $x \in X$ and $y \in Y$. Then we have $g \in \Lip_1(X)$ and $h \in \Lip_1(Y)$.
\end{lem}




\begin{prop}\label{LO}
Let $X$ and $Y$ be two mm-spaces and let $F, G \in \F^2$. If $F \leq G$, that is, $F(s, t) \leq G(s, t)$ for every $s, t$, then we have
\begin{equation}\label{LOeq}
\OD(X \times_F Y ; - \kappa) \leq \OD(X \times_G Y ; -\kappa)
\end{equation}
for any $\kappa > 0$.
\end{prop}

\begin{proof}
We take any $f \in \Lip_1(X \times_F Y)$. Since
\begin{equation*}
|f(x, y) - f(x', y')| \leq F(d_X(x, x'), d_Y(y, y')) \leq G(d_X(x, x'), d_Y(y, y')),
\end{equation*}
we have $f \in \Lip_1(X \times_G Y)$. Thus
\begin{equation*}
\PD(f_*(m_X \otimes m_Y) ; 1-\kappa) \leq \OD(X \times_G Y ; -\kappa),
\end{equation*}
which implies \eqref{LOeq}. The proof is completed.
\end{proof}

\begin{rem}
It is well-known that the observable diameter is monotone with respect to the Lipschitz order which is a partial order relation over $\mathcal{X}$ (see \cite{MMG}*{Proposition 2.18}). Proposition \ref{LO} is a special case of this.
\end{rem}

The following lemmas are keys to the proof of Theorem \ref{main}.

\begin{lem}\label{key_lp}
Let $X$ and $Y$ be two mm-spaces and let $p \in [1, +\infty]$. Then we have
\begin{equation}
\begin{split}
& \OD(X \times_p Y ; - (\kappa + \kappa')) \\
& \leq \OD(X ; - \kappa) + 2\OD(Y ; -\kappa')
\end{split}
\end{equation}
for any $\kappa \in (0, 1)$ and any $\kappa' \in (0, 1/2)$.
\end{lem}

\begin{lem}\label{key_F}
Let $X$ and $Y$ be two mm-spaces and let $F \in \F^2$. Then we have
\begin{equation}
\begin{split}
& \OD(X \times_F Y ; - 2(\kappa + \kappa')) \\
& \leq 4F(\OD(X ; - \kappa), 0) + 8F(0, \OD(Y ; -\kappa'))
\end{split}
\end{equation}
for any $\kappa \in (0, 1)$ and any $\kappa' \in (0, 1/4)$.
\end{lem}

\begin{proof}[Proof of Lemma \ref{key_lp}]
We take any real numbers $\varepsilon, \delta > 0$ with
\begin{equation*}
\OD(X; - \kappa) < \varepsilon \text{ and } \OD(Y; -\kappa') < \delta,
\end{equation*}
and take any $f \in \Lip_1(X \times_p Y)$. It suffices to prove
\begin{equation}\label{goal_eq}
\PD(f_*(m_X \otimes m_Y) ; 1-(\kappa + \kappa')) \leq \varepsilon + 2\delta. 
\end{equation}

For all $x \in X$, let $f_x := f(x, \cdot) \in \Lip_1(Y)$. We define $g \colon X \to \R$ by 
\begin{equation*}
g(x) := \lm(f_x ; m_Y)
\end{equation*}
for $x \in X$. By Lemma \ref{lm_lip}, we have $g \in \Lip_1(X)$. By $\OD(X; - \kappa) < \varepsilon$, we have $\PD(g_* m_X ; 1- \kappa) < \varepsilon$. Thus there exists an interval $I \subset \R$ such that 
\begin{equation*}
g_* m_X(I) \geq 1-\kappa \text{ and } \diam{I} < \varepsilon.
\end{equation*}
We take countable dense points $\{x_i \}_{i \in \N}$ in $g^{-1}(I)$. By $\OD(Y; -\kappa') < \delta$, we have $\PD((f_{x_i})_* m_Y ; 1- \kappa') < \delta$ for any $i \in \N$. For each $i \in \N$, there exists an interval $J_i \subset \R$ such that 
\begin{equation*}
(f_{x_i})_* m_Y(J_i) \geq 1-\kappa' \text{ and } \diam{J_i} < \delta. 
\end{equation*}
Combining $(f_{x_i})_* m_Y(J_i) > 1/2$ and the connectedness of $J_i$ leads to $\lm(f_{x_i} ; m_Y) \in J_i$. 

We take any $\eta > 0$ and fix it. We define a Borel set $\tilde{A}_\eta$ by
\begin{equation}
\tilde{A}_\eta := \bigsqcup_{i \in \N} \left(U_\eta^i \times f_{x_i}^{-1}(J_i) \right) \subset X \times Y,
\end{equation}
where $U_\eta^1 := U_\eta(x_1)$ and $U_\eta^i := U_\eta(x_i) \setminus \bigcup_{j < i} U_\eta^j$ for $i \geq 2$. Then we see that
\begin{align*}
m_X \otimes m_Y(\tilde{A}_\eta) & = \sum_{i = 1}^\infty m_X(U_\eta^i) m_Y(f_{x_i}^{-1}(J_i)) \geq (1 - \kappa') \sum_{i = 1}^\infty m_X(U_\eta^i) \\
& \geq (1 - \kappa') m_X(g^{-1}(I)) \geq (1-\kappa)(1 - \kappa') \geq 1-(\kappa + \kappa').
\end{align*}
A closed subset $A_\eta$ of $\R$ is defined as the closure of $f(\tilde{A}_\eta)$. Then we have
\begin{equation*}
f_*(m_X \otimes m_Y)(A_\eta) \geq m_X \otimes m_Y(\tilde{A}_\eta) \geq 1-(\kappa + \kappa').
\end{equation*}
Moreover, we have
\begin{equation*}
\diam{A_\eta} = \diam{f(\tilde{A}_\eta)} \leq \varepsilon + 2\delta + 2\eta.
\end{equation*}
In fact, for any $(x, y) \in U_\eta^i \times f_{x_i}^{-1}(J_i)$ and $(x', y') \in U_\eta^j \times f_{x_j}^{-1}(J_j)$,
\begin{align*}
& \, |f(x,y) - f(x', y')| \\
\leq & \, |f(x,y) - f(x_i, y)| + |f(x_i,y) - g(x_i)| + |g(x_i) - g(x_j)| \\
& + |g(x_j) - f(x_j, y')| +|f(x_j,y') - f(x', y')| \\
\leq & \, d_X(x, x_i) + \diam{J_i} + \diam{I} + \diam{J_j} + d_X(x_j, x') \\
\leq & \, \varepsilon + 2\delta + 2\eta.
\end{align*}
Thus we have $\PD(f_*(m_X \otimes m_Y) ; 1-(\kappa + \kappa')) \leq \diam{A_\eta} \leq \varepsilon + 2\delta+ 2\eta$, which implies (\ref{goal_eq}). This completes the proof.
\end{proof}

\begin{proof}[Proof of Lemma \ref{key_F}]
Let $\kappa$ and $\kappa'$ be two real numbers with $\kappa \in (0, 1)$ and $\kappa' \in (0, 1/4)$. Let $F \in \F^2$. The two functions $F_1$ and $F_2$ on $[0, +\infty)$ defined by
\begin{equation*}
F_1(s) := F(s, 0), \quad F_2(t) := F(0, t)
\end{equation*}
for $s, t \in [0, +\infty)$ belong to $\F^1$, and then the function $G$ on $[0, +\infty)^2$ defined by
\begin{equation*}
G(s, t) := F_1(s) + F_2(t) =F(s, 0) + F(0, t)
\end{equation*}
for $s, t \in [0, +\infty)$ belongs to $\F^2$. Since $F \leq G$ holds, by Proposition \ref{LO}, we have
\begin{equation*}
\OD(X \times_F Y ; - 2(\kappa + \kappa')) \leq \OD(X \times_G Y ; -2(\kappa + \kappa')).
\end{equation*}
Taking into account that $X \times_G Y$ is mm-isomorphic to 
\begin{equation*}
(X, F_1 \circ d_X, m_X) \times_1 (Y, F_2 \circ d_Y, m_Y), 
\end{equation*}
by Lemma \ref{key_lp} and Lemma \ref{key_1dim} (2), we have
\begin{equation*}
\begin{split}
& \OD(X \times_G Y ; - 2(\kappa + \kappa')) \\
& \leq 4F_1(\OD(X ; -\kappa)) + 8F_2(\OD(Y ; -\kappa')) \\
& = 4F(\OD(X ; -\kappa), 0) + 8F(0, \OD(Y ; -\kappa')).
\end{split}
\end{equation*}
This completes the proof.
\end{proof}

\begin{ex}\label{sphere_norm}
Let $S^n(r_n)$, $n = 1, 2, \ldots$, be the sphere of radius $r_n > 0$ in $\R^{n+1}$ and let $\sigma^n$ be the Riemannian volume measure on $S^n(r_n)$ normalized as $\sigma^n(S^n(r_n)) = 1$. Let $F_n$ be the function in $\F^1$ defined by
\begin{equation*}
F_n(s) := \left\{\begin{array}{ll} 2r_n \sin{\frac{s}{2r_n}} & \text{if } s \leq \pi r_n, \\ 2r_n & \text{if } s > \pi r_n. \end{array}\right.
\end{equation*}
for $s \in [0, +\infty)$. Let $d$ be the Riemannian distance on $S^n(r_n)$ and let $\| \cdot \|$ be the Euclidean norm on $\R^{n+1}$. We see that 
\begin{equation*}
\| x - x' \|  = F_n(d(x, x')) \leq d(x, x')
\end{equation*}
for any $x, x' \in S^n(r_n)$. Thus, by Proposition \ref{LO}, we have
\begin{equation*}
\OD((S^n(r_n), \| \cdot \|, \sigma^n); -\kappa) \leq \OD((S^n(r_n), d, \sigma^n); -\kappa)
\end{equation*}
for any $\kappa > 0$. The sequence $\{(S^n(r_n), \| \cdot \|, \sigma^n)\}_{n \in \N}$ is a L\'evy family if $r_n = o(\sqrt{n})$. We use this example in Section 5.
\end{ex}

\section{Convergence of product spaces}

\subsection{Box-convergence of product spaces}

The purpose of this subsection is to prove the following two propositions and their corollaries. We denote by $\N$ the set of positive integers and by $\F^N$ the set of continuous metric preserving functions on $[0, +\infty)^N$.

\begin{prop}\label{box1}
Let $X$, $Y$, $Z$, and $W$ be four mm-spaces and let $F \in \F^2$. Then we have
\begin{equation}\label{box1_eq1}
\begin{split}
& \square(X \times_F Z, Y \times_F W)  \\
& \leq \max{\left\{ \square(X, Y) + \square(Z, W), 2F(\frac{1}{2}\square(X, Y), \frac{1}{2}\square(Z, W)) \right\}}.
\end{split}
\end{equation}
In particular, for any $p \in [1, +\infty]$, it holds that
\begin{equation}\label{box1_eq2}
\square(X \times_p Z, Y \times_p W) \leq \square(X, Y) + \square(Z, W).
\end{equation}
\end{prop}

\begin{prop}\label{box2}
Let $\{X_n\}_{n \in \N}$ and $\{Y_n\}_{n \in \N}$ be two sequences of mm-spaces $\square$-converging to mm-spaces $X$ and $Y$ respectively. Let $F_n$, $n = 1, 2, \ldots$, and $F$ be functions in $\F^2$ such that $F_n$ converges pointwise to $F$ as $n \to \infty$. Then the sequence $\{X_n \times_{F_n} Y_n\}_{n \in \N}$ of their product spaces $\square$-converges to the product space $X \times_F Y$.
\end{prop}

\begin{cor}\label{box3}
Let $\{X_n\}_{n \in \N}$ and $\{Y_n\}_{n \in \N}$ be two sequences of mm-spaces $\square$-converging to mm-spaces $X$ and $Y$ respectively. Assume that $p_n \in [1, +\infty]$ converges to $p \in [1, +\infty]$ as $n \to \infty$. Then the sequence $\{X_n \times_{p_n} Y_n\}_{n \in \N}$ of their $l_{p_n}$-product spaces $\square$-converges to the $l_p$-product space $X \times_p Y$.
\end{cor}

\begin{cor}\label{box4}
Let $\{X_n\}_{n \in \N}$ be a sequence of mm-spaces $\square$-converging to an mm-space $X$. Let $F_n$, $n=1, 2, \ldots$, and $F$ be functions in $\F^1$ such that $F_n$ converges pointwise to $F$ as $n \to \infty$. Then the sequence $\{(X_n, F_n \circ d_{X_n}, m_{X_n})\}_{n \in \N}$ of mm-spaces $\square$-converges to the mm-space $(X, F \circ d_X, m_X)$.
\end{cor}

Given two maps $f \colon X \to X'$ and $g \colon Y \to Y'$, a map $f \times g \colon X \times Y \to X' \times Y'$ is defined by 
\begin{equation}\label{timesmap}
(f \times g)(x, y) := (f(x), g(y))
\end{equation}
for $(x, y) \in X \times Y$. In this paper, for two maps $f$ and $g$, we always use the notation $f \times g$ in the sense of (\ref{timesmap}).

\begin{proof}[Proof of Proposition \ref{box1}]
This proposition is a generalization of \cite{OS}*{Lemma 3.3} and its proof is similar to that of \cite{OS}*{Lemma 3.3}.

Take any two real numbers $\varepsilon, \delta > 0$ such that $\square(X, Y) < \varepsilon$, $\square(Z, W) < \delta$. Then there exist four parameters $\varphi \colon I \to X$, $\psi \colon I \to Y$, $\xi \colon I \to Z$, and $\eta \colon I \to W$ of every $X$, $Y$, $Z$, and $W$ respectively, and two Borel subsets $I_1$, $I_2 \subset I$ such that
\begin{align*}
& |d_X(\varphi(s), \varphi(t)) - d_Y(\psi(s), \psi(t))| < \varepsilon \text{ for any } s,t \in I_1, \\
& |d_Z(\xi(s), \xi(t)) - d_W(\eta(s), \eta(t))| < \delta \text{ for any } s,t \in I_2, \\
& \mathcal{L}^1(I_1) \geq 1- \varepsilon, \text{ and } \mathcal{L}^1(I_2) \geq 1- \delta.
\end{align*}
Let $\Phi$ be a parameter of $(I \times I, \mathcal{L}^1 \otimes \mathcal{L}^1)$, i.e., a Borel measurable map $\Phi \colon I \to I \times I$ with
\begin{equation*}
\Phi_* \mathcal{L}^1 = \mathcal{L}^1 \otimes \mathcal{L}^1,
\end{equation*}
where the existence of $\Phi$ follows from Lemma \ref{prm_exist}. Then the Borel measurable map $(\varphi \times \xi) \circ \Phi \colon I \to X \times Z$ is a parameter of $X \times Z$. In fact, we have
\begin{equation*}
\left((\varphi \times \xi) \circ \Phi \right)_* \mathcal{L}^1 = (\varphi \times \xi)_* (\mathcal{L}^1 \otimes \mathcal{L}^1) = \varphi_*\mathcal{L}^1 \otimes \xi_* \mathcal{L}^1 = m_X \otimes m_Z.
\end{equation*}
Similarly, the Borel measurable map $(\psi \times \eta) \circ \Phi \colon I \to Y \times W$ is a parameter of $Y \times W$. Setting $I_0 := \Phi^{-1}(I_1 \times I_2)$, we have 
\begin{equation*}
\mathcal{L}^1(I_0) = \mathcal{L}^1(\Phi^{-1}(I_1 \times I_2)) = \mathcal{L}^1(I_1)\mathcal{L}^1(I_2) \geq (1 - \varepsilon)(1 - \delta) \geq 1 - (\varepsilon + \delta).
\end{equation*}
In addition, we define four functions $\varphi', \psi', \xi', \eta'$ by
\begin{equation*}
(\varphi'(s), \xi'(s)) = (\varphi \times \xi) \circ \Phi(s) \text{ and } (\psi'(s), \eta'(s)) = (\psi \times \eta) \circ \Phi(s)
\end{equation*}
for any $s \in I$. Then, for any $s, t \in I_0$, we have
\begin{align*}
&|d_F((\varphi \times \xi) \circ \Phi(s), (\varphi \times \xi) \circ \Phi(t)) - d_F((\psi \times \eta) \circ \Phi(s), (\psi \times \eta) \circ \Phi(t))| \\
& = |F(d_X(\varphi'(s), \varphi'(t)), d_Z(\xi'(s), \xi'(t))) - F(d_Y(\psi'(s), \psi'(t)), d_W(\eta'(s), \eta'(t)))| \\
& \leq F(|d_X(\varphi'(s), \varphi'(t)) - d_Y(\psi'(s), \psi'(t))|, |d_Z(\xi'(s), \xi'(t)) - d_W(\eta'(s), \eta'(t))|) \\
& \leq 2F(\frac{1}{2}\varepsilon, \frac{1}{2}\delta),
\end{align*}
where the inequalities follow from Corollary \ref{cor_mpf}. Thus we have
\begin{equation}\label{box1pf_eq}
\square(X \times_F Z, Y \times_F W) \leq \max{\left\{ \varepsilon + \delta, 2F(\frac{1}{2}\varepsilon, \frac{1}{2}\delta) \right\}},
\end{equation}
so that we obtain (\ref{box1_eq1}). Furthermore, in the case of $F = F_p$ in Example \ref{exF} (1) for $p \in [1, +\infty]$, we see that $F_p(s, t) \leq s + t$ for every $s, t$. Combining this and (\ref{box1pf_eq}) implies 
\begin{equation*}
\square(X \times_F Z, Y \times_F W) \leq \varepsilon + \delta,
\end{equation*}
which means (\ref{box1_eq2}). The proof of the proposition is completed.
\end{proof}

We need the following two lemmas for the proof of Proposition \ref{box2}.

\begin{lem}\label{lprok}
Let $\mu, \mu'$ be two Borel probability measures on a separable metric space $X$ and let $\nu, \nu'$ be two Borel probability measures on a separable metric space $Y$. Let $F$ be a function in $\F^2$. Then we have
\begin{equation}
\begin{split}
&\lprok{\lambda}(\mu \otimes \nu, \mu' \otimes \nu') \\
&\leq  \max{\left\{ \lprok{\lambda}(\mu, \mu') + \lprok{\lambda}(\nu, \nu'), 2F(\lprok{\lambda}(\mu, \mu'), \lprok{\lambda}(\nu, \nu')) \right\}} \label{lprok_eq}
\end{split}
\end{equation}
for any $\lambda > 0$, where $\lprok{\lambda}$ in the left-hand side is with respect to $d_F$.
\end{lem}

\begin{proof}
We take any two real numbers $\varepsilon$ and $\delta$ such that $\lprok{\lambda}(\mu, \mu') < \varepsilon$, $\lprok{\lambda}(\nu, \nu') < \delta$ and fix them. For the proof of (\ref{lprok_eq}), it suffices to prove that
\begin{equation}\label{lproklemgoal}
(\mu \otimes \nu)(A) \leq (\mu' \otimes \nu') (U_{2F(\varepsilon, \delta)+\eta}(A)) + \lambda (\varepsilon + \delta) 
\end{equation}
for any Borel subset $A \subset X \times Y$ and for any $\eta > 0$. The famous $\pi$-$\lambda$ theorem shows that if we prove
\begin{equation}\label{lproklemgoal2}
(\mu \otimes \nu)(B \times C) \leq (\mu' \otimes \nu') (U_{2F(\varepsilon, \delta)+\eta}(B \times C)) + \lambda (\varepsilon + \delta) 
\end{equation}
for any Borel subsets $B \subset X$ and $C \subset Y$, then we obtain (\ref{lproklemgoal}) for any Borel subset $A \subset X \times Y$. Let us prove (\ref{lproklemgoal2}). By $\lprok{\lambda}(\mu, \mu') < \varepsilon$ and $\lprok{\lambda}(\nu, \nu') < \delta$, we have
\begin{align*}
& (\mu \otimes \nu)(B \times C) \left(= \mu(B) \nu(C) \right) \\
& \leq \left(\mu'(U_{\varepsilon}(B)) + \lambda\varepsilon\right) \nu(C) \\
& \leq \mu'(U_{\varepsilon}(B))\nu(C)  + \lambda\varepsilon \\ 
& \leq \mu'(U_{\varepsilon}(B)) \left(\nu'(U_{\delta}(C)) + \lambda\delta \right)  + \lambda\varepsilon \\
& \leq \mu'(U_{\varepsilon}(B)) \nu'(U_{\delta}(C)) + \lambda(\varepsilon + \delta) \\
& \leq (\mu' \otimes \nu') (U_{2F(\varepsilon, \delta)+\eta}(B \times C)) + \lambda (\varepsilon + \delta) 
\end{align*}
for any Borel subsets  $B \subset X$ and $C \subset Y$, where the last inequality follows from Corollary \ref{cor_mpf} (2). Thus we obtain (\ref{lproklemgoal2}) and then (\ref{lprok_eq}). The proof is completed.
\end{proof}

\begin{lem}\label{ptcpt}
Let $F_n$, $n = 1, 2, \ldots$, and $F$ be metric preserving functions. If $F$ is continuous and $F_n$ converges pointwise to $F$, then $F_n$ uniformly converges to $F$ on compact sets.
\end{lem}

\begin{proof}
We take any compact set $K \subset [0, +\infty)^2$ and any real number $\varepsilon > 0$. Let us prove that
\begin{equation}\label{ptcpt_eq}
\sup_{(s, t) \in K}|F_n(s, t) - F(s, t)| \leq 7\varepsilon
\end{equation}
holds for every sufficiently large $n$. By the continuity of $F$, there exists a real number $\delta > 0$ such that $F(\delta, \delta) \leq \varepsilon$. By the compactness of $K$, we find finite points $\{(s_i, t_i)\}_{i=1}^k$ in $K$ such that
\begin{equation*}
K \subset \bigcup^k_{i=1} \left( U_\delta(s_i) \times U_\delta(t_i) \right).
\end{equation*}
Let $N \in \N$ be a number such that
\begin{equation*}
\max_{i=1, \ldots, k} |F_n(s_i, t_i) - F(s_i, t_i)| \leq \varepsilon \text{ and } |F_n(\delta, \delta) - F(\delta, \delta)| \leq \varepsilon
\end{equation*}
hold for all $n \geq N$. Given a fixed point $(s, t) \in K$, we find $i \in \{1, \ldots, k\}$ such that $(s, t) \in U_\delta(s_i) \times U_\delta(t_i)$. By Corollary \ref{cor_mpf}, we have
\begin{align*}
&|F_n(s, t) - F(s, t)| \\
&\leq |F_n(s_i, t_i) - F(s_i, t_i)| + F_n(|s - s_i|, |t - t_i|) + F(|s - s_i|, |t - t_i|) \\
&\leq |F_n(s_i, t_i) - F(s_i, t_i)| + 2F_n(\delta, \delta) + 2F(\delta, \delta) \\
&\leq |F_n(s_i, t_i) - F(s_i, t_i)| + 2|F_n(\delta, \delta) - F(\delta, \delta)| + 4F(\delta, \delta) \\
&\leq 7\varepsilon
\end{align*}
for every $n \geq N$. Thus we obtain \eqref{ptcpt_eq}. This completes the proof.
\end{proof}

\begin{proof}[Proof of Proposition \ref{box2}]
It is sufficient to prove that, for any real number $\varepsilon > 0$, there exists $N(\varepsilon) \in \N$ such that an $\varepsilon$-mm-isomorphism $\Phi_n \colon X_n \times_{F_n} Y_n \to X \times_F Y$ exists for each $n \geq N(\varepsilon)$. Take any $\varepsilon > 0$ and fix it. By the inner regularity of probability measures $m_X$ and $m_Y$, there exist compact sets $K \subset X$ and $K' \subset Y$ such that
\begin{equation*}
m_X(K) \geq 1- \varepsilon \text{ and } m_Y(K') \geq 1- \varepsilon.
\end{equation*}
Let $D_\varepsilon := \max{\{ \diam{K}, \diam{K'}\}} + 3\varepsilon$. Since $\{X_n\}_{n \in \N}$, $\{Y_n\}_{n \in \N}$ $\square$-converge to $X$, $Y$ respectively and $F_n$ uniformly converges to $F$ on $[0, D_\varepsilon]^2$ following from Lemma \ref{ptcpt}, there exists $N(\varepsilon) \in \N$ such that, for any $n \geq N(\varepsilon)$,
\begin{itemize}
\item an $\varepsilon$-mm-isomorphism $f_n \colon X_n \to X$ exists,
\item an $\varepsilon$-mm-isomorphism $g_n \colon Y_n \to Y$ exists,
\item $|F_n(s, t) - F(s, t)| < \varepsilon$ for all $s, t \in [0, D_\varepsilon]$.
\end{itemize}
If we prove that the map $f_n \times g_n$ is an $(4F(\varepsilon, \varepsilon) + 6\varepsilon)$-mm-isomorphism for any $n \geq N(\varepsilon)$, then we obtain the conclusion with $\Phi_n = f_n \times g_n$. Take any $n \geq N(\varepsilon)$ and fix it. Let us prove that the map $f_n \times g_n$ is an $(4F(\varepsilon, \varepsilon) + 6\varepsilon)$-mm-isomorphism. By Lemma \ref{lprok}, taking 
\begin{equation*}
(f_n \times g_n)_* (m_{X_n} \otimes m_{Y_n}) = {f_n}_* m_{X_n} \otimes {g_n}_* m_{Y_n}
\end{equation*}
into account, we have
\begin{equation}\label{proof2_eq1}
\prok((f_n \times g_n)_* (m_{X_n} \otimes m_{Y_n}), m_X \otimes m_Y) \leq \max{\{ 2\varepsilon, 4F(\varepsilon, \varepsilon)\}}.
\end{equation}
Let $X'_n$, $Y'_n$ be nonexceptional domains of $f_n$, $g_n$ respectively and let
\begin{equation*}
\widetilde{X}_n := X'_n \cap f_n^{-1}(U_\varepsilon(K)) \text{ and } \widetilde{Y}_n := Y'_n \cap g_n^{-1}(U_\varepsilon(K')).
\end{equation*}
Since $\prok({f_n}_* m_{X_n}, m_X) \leq \varepsilon$, we have
\begin{align*}
m_{X_n}(\widetilde{X}_n) & \geq m_{X_n}(X'_n) + m_{X_n}(f_n^{-1}(U_\varepsilon(K))) - 1 \\
& \geq m_{X_n}(X'_n) + ( m_X(K) - \varepsilon ) - 1 \\
& \geq 1 - 3\varepsilon.
\end{align*}
Similarly, $\prok({g_n}_* m_{Y_n}, m_Y) \leq \varepsilon$ implies $m_{Y_n}(\widetilde{Y}_n) \geq 1 - 3\varepsilon$. Thus we obtain
\begin{equation}\label{proof2_eq2}
m_{X_n} \otimes m_{Y_n}(\widetilde{X}_n \times \widetilde{Y}_n) \geq 1 - 6\varepsilon.
\end{equation}
Furthermore we see that $\diam{\widetilde{X}_n} \leq D_\varepsilon$ and $\diam{\widetilde{Y}_n} \leq D_\varepsilon$. In fact,
\begin{equation*}
d_{X_n}(x, x') \leq d_X(f_n(x), f_n(x')) + \varepsilon \leq \diam{K} + 3\varepsilon \leq D_\varepsilon
\end{equation*}
for any $x, x' \in \widetilde{X}_n$. It is possible to check $\diam{\widetilde{Y}_n} \leq D_\varepsilon$ similarly. Therefore, for any $(x, y), (x', y') \in \widetilde{X}_n \times \widetilde{Y}_n$, we have 
\begin{align*}
& \,|d_{F_n}((x, y), (x', y')) - d_F((f_n \times g_n)(x, y), (f_n \times g_n)(x', y'))| \\
= & \, |F_n(d_{X_n}(x, x'), d_{Y_n}(y, y')) - F(d_X(f_n(x), f_n(x')), d_Y(g_n(y), g_n(y')))| \\
\leq & \, |F_n(d_{X_n}(x, x'), d_{Y_n}(y, y')) - F(d_{X_n}(x, x'), d_{Y_n}(y, y'))| \\
& + F(|d_{X_n}(x, x') - d_X(f_n(x), f_n(x'))|, |d_{Y_n}(y, y') - d_Y(g_n(y), g_n(y'))|) \\
\leq & \, \varepsilon + 2F(\varepsilon, \varepsilon).
\end{align*}
Combining this with (\ref{proof2_eq1}) and (\ref{proof2_eq2}) means that the map $f_n \times g_n$ is an $(4F(\varepsilon, \varepsilon) + 6\varepsilon)$-mm-isomorphism. The proof of the proposition is completed.
\end{proof}

\begin{proof}[Proof of Corollary \ref{box3}]
We apply Proposition \ref{box2} with $F_n = F_{p_n}$, $n = 1, 2,\ldots$, and $F = F_p$, where $F_p$ is the function of Example \ref{exF} (\ref{exF1}). 
\end{proof}

\begin{proof}[Proof of Corollary \ref{box4}]
Let $Y_n$, $n = 1, 2, \ldots$, and $Y$ be one-point mm-spaces and let $G_n$ and $G$ be the functions on $[0, +\infty)^2$ defined by
\begin{equation*}
G_n(s, t) := F_n(s) + t \text{ and } G(s, t) := F(s) + t
\end{equation*}
for $s, t \in [0, + \infty)$. We just apply Proposition \ref{box2}. Note that the mm-space $X_n \times_{G_n} Y_n$ is mm-isomorphic to $(X_n, F_n \circ d_{X_n}, m_{X_n})$.
\end{proof}

\subsection{Concentration of product spaces}
Our goals in this subsection are to prove the half of Theorem \ref{main}, and to obtain Corollary \ref{main_cor} and the half of Theorem \ref{main_1dim} as its corollaries.

\begin{dfn}[$1$-Lipschitz up to an additive error]
Let $X$ be an mm-space and $Y$ be a metric space. A map $f \colon X \to Y$ is said to be 1-{\it Lipschitz up to} ({\it an additive error}) {\it $\varepsilon \geq 0$} if there exists a Borel subset $X_0 \subset X$ such that
\begin{enumerate}
\item $m_X(X_0) \geq 1 - \varepsilon$,
\item $d_Y(f(x), f(x')) \leq d_X(x, x') + \varepsilon$ for any $x, x' \in X_0$.
\end{enumerate}
We call such a set $X_0$ a {\it nonexceptional domain} of $f$.
\end{dfn}

\begin{lem}[\cite{MMG}*{Lemma 5.4}]\label{MHext}
If a function $f \colon X \to \R$ on an mm-space $X$ is $1$-Lipschitz up to an additive error $\varepsilon \geq 0$, then there exists a $1$-Lipschitz function $\tilde{f} \colon X \to \R$ such that
\begin{equation*}
\kf^{m_X}(f, \tilde{f}) \leq \varepsilon.
\end{equation*}
\end{lem}

\begin{lem}[\cite{MMG}*{Lemma 5.27}]\label{mmg5.27}
Let $X$ and $Y$ be two mm-spaces and $p \colon X \to Y$ a Borel measurable map. For two real numbers $\varepsilon, \delta > 0$, we consider the two following conditions.
\begin{enumerate}
\renewcommand{\labelenumi}{{\rm(A$_\varepsilon$)}}
\item  $p^* \Lip_1(Y) \subset U_{\varepsilon + \eta}(\Lip_1(X))$ for any $\eta > 0$.
\renewcommand{\labelenumi}{{\rm(B$_\delta$)}}
\item  $p$ is $1$-Lipschitz up to $\delta$.
\end{enumerate}
Then we have the following {\rm (1)} and {\rm (2)}.
\begin{enumerate}
\item There exists a real number $\delta = \delta(Y, \varepsilon) > 0$ for any $\varepsilon > 0$ such that $\lim_{\varepsilon \to 0}\delta(Y, \varepsilon) = 0$ and if {\rm (A$_\varepsilon$)} holds and if $\prok(p_* m_X, m_Y) < \varepsilon$, then we have {\rm (B$_\delta$)}.
\item If {\rm (B$_\delta$)} holds, then we have {\rm (A$_\delta$)}.
\end{enumerate}
\end{lem}

The following lemma gives an condition equivalent to (\ref{main2}) of Theorem \ref{main}. The condition (\ref{main2}) of Theorem \ref{main} strengthen seemingly.

\begin{lem}\label{equiv_main}
Let $F_n \colon [0, +\infty)^2 \to [0, +\infty)$ be a function, $n = 1, 2, \ldots$. Assume that $F_n$ uniformly converges to a continuous function $F$ on compact sets. Then the following {\rm (1)} and {\rm (2)} are equivalent to each other.
\begin{enumerate}
\item For any $s, t \in [0, + \infty)$, 
\begin{equation*}
\lim_{n \to \infty} (F_n(s, t) - \inf_{s\leq s' \text{and } t \leq t'} F_n(s', t')) = 0.
\end{equation*}
\item For any $D > 0$,
\begin{equation*}
\lim_{n \to \infty} \sup_{0 \leq s,t \leq D} (F_n(s, t) - \inf_{s\leq s' \text{and } t \leq t'} F_n(s', t')) = 0.
\end{equation*}
\end{enumerate}
\end{lem}

\begin{proof}
It is trivial that (2) implies (1). We prove that (1) implies (2). Suppose that the condition (2) does not hold in order to prove the contraposition. There exists a real number $D > 0$ such that
\begin{equation*}
\limsup_{n \to \infty} \sup_{0 \leq s,t \leq D} (F_n(s, t) - \inf_{s\leq s' \text{and } t \leq t'} F_n(s', t')) > 0.
\end{equation*}
Choosing a subsequence of $n$, we can assume that there exist a real number $\eta > 0$ and a sequence $\{(s_n, t_n)\}_{n \in \N} \subset [0 ,D]^2$ such that
\begin{equation*}
F_n(s_n, t_n) - \inf_{s_n\leq s' \text{and } t_n \leq t'} F_n(s', t') > \eta.
\end{equation*}
Choosing a subsequence again, we can assume that $s_n$, $t_n$ converge to $s_\infty$, $t_\infty$, respectively, as $n \to \infty$. We see that $s_\infty, t_\infty \leq D$. By the continuity of $F$, there exists a real number $\delta > 0$ such that
\begin{equation*}
|F(s, t) - F(s_\infty, t_\infty)| < \frac{\eta}{8}
\end{equation*}
for any $s, t \in [0, +\infty)$ with $|s-s_\infty|, |t-t_\infty| \leq \delta$. Since $F_n$ uniformly converges to $F$ on $[0, D+\delta]^2$, for every sufficiently large $n$,
\begin{equation*}
\sup_{0 \leq s, t \leq D +\delta} |F_n(s, t) - F(s, t)| < \frac{\eta}{8}.
\end{equation*}
Let $\Gamma := \left\{(s, t) \in  [0, +\infty)^2 \midd |s-s_\infty|, |t-t_\infty| \leq \delta \right\}$. For every sufficiently large $n$ and for every $(s, t), (s', t') \in \Gamma$, we have
\begin{align*}
& \, |F_n(s, t) - F_n(s', t')| \\
\leq & \, |F_n(s, t) - F(s, t)| + |F(s, t) - F(s_\infty, t_\infty)| \\
& + |F(s_\infty, t_\infty) - F(s', t')| + |F(s', t') - F_n(s', t')| \\
< & \, \frac{\eta}{2}.
\end{align*}
Let $s_* := \max\{s_\infty - \delta, 0\}$ and $t_* := \max\{t_\infty - \delta, 0\}$. Taking into account that $(s_*, t_*) \in \Gamma$ and $(s_n, t_n) \in \Gamma$ for every sufficiently large $n$, we have
\begin{align*}
& F_n(s_*, t_*) - \inf_{s_* \leq s' \text{and } t_* \leq t'} F_n(s', t') \\
& > F_n(s_n, t_n) - \frac{\eta}{2} - \inf_{s_n \leq s' \text{and } t_n \leq t'} F_n(s', t') > \frac{\eta}{2}
\end{align*}
for every sufficiently large $n$. This means that the condition (1) does not hold. The proof is completed.
\end{proof}

\begin{lem}\label{lip_up}
Let $X_n$, $Y_n$, $X$, and $Y$ be mm-spaces and let $F_n$ and $F$ be functions in $\F^2$, where $n = 1, 2, \ldots$. Assume that $F_n$ converges pointwise to $F$ and satisfies the condition {\rm (\ref{main2})} of Theorem \ref{main}. Let $p_n \colon X_n \to X$ and $q_n \colon Y_n \to Y$ be maps. If $p_n$, $q_n$ are $1$-Lipschitz up to $\varepsilon_n$, $\delta_n$ respectively, and both 
\begin{equation*}
\prok({p_n}_* m_{X_n}, m_X) \leq \varepsilon_n \text{ and } \prok({q_n}_* m_{Y_n}, m_Y) \leq \delta_n
\end{equation*}
hold for some sequences $\varepsilon_n, \delta_n \to 0$ as $n \to \infty$, then the map $p_n \times q_n \colon X_n \times_{F_n} Y_n \to X \times_F Y$ is $1$-Lipschitz up to $\eta_n$ for some sequence $\eta_n \to 0$ as $n \to \infty$.
\end{lem}

\begin{proof}
The proof is similar to that of Proposition \ref{box2}.

Take any real number $\varepsilon > 0$. It suffices to prove that the map $p_n \times q_n$ is $1$-Lipschitz up to $F(\varepsilon, \varepsilon) + 6\varepsilon$ for any sufficiently large $n \in \N$. By the inner regularity of $m_X$ and $m_Y$, there exist compact sets $K \subset X$ and $K' \subset Y$ such that
\begin{equation*}
m_X(K) \geq 1- \varepsilon \text{ and } m_Y(K') \geq 1- \varepsilon.
\end{equation*}
Let $D_\varepsilon := \max{\{ \diam{K}, \diam{K'}\}} + 2\varepsilon$. Then, by the assumptions and Lemma \ref{equiv_main}, for any sufficiently large $n \in \N$,
\begin{itemize}
\item $p_n, q_n$ are both $1$-Lipschitz up to $\varepsilon$,
\item $\prok({p_n}_* m_{X_n}, m_X) \leq \varepsilon$ and $\prok({q_n}_* m_{Y_n}, m_Y) \leq \varepsilon$ hold,
\item $|F_n(s, t) - F(s, t)| < \varepsilon$ holds for all $s, t \in [0, D_\varepsilon]$,
\item $F_n(s, t) \leq F_n(s', t') + \varepsilon$ holds for any $s, t \in [0, D_\varepsilon]$ and for any $s', t' \in [0, +\infty)$ with $s \leq s'$ and $t \leq t'$.
\end{itemize}
Let $X'_n$, $Y'_n$ be nonexceptional domains of $p_n$, $q_n$ respectively and let
\begin{equation*}
\widetilde{X}_n := X'_n \cap p_n^{-1}(U_\varepsilon(K)) \text{ and } \widetilde{Y}_n := Y'_n \cap q_n^{-1}(U_\varepsilon(K')).
\end{equation*}
By the similar proof to that of (\ref{proof2_eq2}), we have 
\begin{equation*}
m_{X_n} \otimes m_{Y_n}(\widetilde{X}_n \times \widetilde{Y}_n) \geq 1 - 6\varepsilon.
\end{equation*}
For any $(x, y), (x', y') \in \widetilde{X}_n \times \widetilde{Y}_n$, we have
\begin{align*}
&d_F((\varphi_n \times \psi_n)(x, y), (\varphi_n \times \psi_n)(x', y')) \\
= & \, F(d_X(\varphi_n(x), \varphi_n(x')), d_Y(\psi_n(y), \psi_n(y'))) \\
\leq & \, F_n(d_X(\varphi_n(x), \varphi_n(x')), d_Y(\psi_n(y), \psi_n(y'))) +\varepsilon \\
\leq & \, F_n(d_{X_n}(x, x') + \varepsilon, d_{Y_n}(y, y') + \varepsilon) + 2\varepsilon \\
\leq & \, F_n(d_{X_n}(x, x'), d_{Y_n}(y, y')) + F_n(\varepsilon, \varepsilon) + 2\varepsilon, \\
\leq & \, d_{F_n}((x, y), (x', y')) + F(\varepsilon, \varepsilon) + 3\varepsilon,
\end{align*}
where the first and second inequalities follow from
\begin{equation*}
d_X(\varphi_n(x), \varphi_n(x')) \leq D_\varepsilon, \quad d_Y(\psi_n(y), \psi_n(y')) \leq D_\varepsilon.
\end{equation*}
Therefore the map $p_n \times q_n$ is $1$-Lipschitz up to $F(\varepsilon,  \varepsilon) + 6\varepsilon$. This completes the proof. 
\end{proof}

Given two subsets $A$ and $B$ of a metric space $X$, we define
\begin{equation*}
d_X(A, B) := \inf_{a \in A, \, b \in B} d_X(a, b).
\end{equation*}

\begin{dfn}[$\kappa$-distance]
Let $\kappa > 0$ and let $X$ be an mm-space. We define the $\kappa$-distance $d_+(A_1, A_2; + \kappa)$ between two Borel subsets $A_1$ and $A_2$ of $X$ as the supremum of $d_X(B_1, B_2)$ over all Borel subsets $B_1 \subset A_1$ and $B_2 \subset A_2$ with $m_X(B_1) \geq \kappa$ and $m_X(B_2) \geq \kappa$. We set $d_+(A_1, A_2; + \kappa) := 0$ if $\min\{m_X(A_1), m_X(A_2)\} < \kappa$.
\end{dfn}

\begin{thm}[Fibration theorem, \cite{Grmv}*{3.$\frac{1}{2}$.47. Proposition}, \cite{MMG}*{Theorem 9.8}]\label{Fib}
Let $p_n \colon X_n \to X$ be a Borel measurable map between mm-spaces $X_n$ and $X$, where $n=1, 2, \ldots$, such that $\prok({p_n}_* m_{X_n}, m_X)$ tends to $0$ as $n \to \infty$. Then, each $p_n$ enforces $\varepsilon_n$-concentration of $X_n$ to $X$ for some sequence $\varepsilon_n \to 0$ if and only if we have the following {\rm (1)}, {\rm (2)}, and {\rm (3)}.
\begin{enumerate}
\item Each $p_n$ is 1-Lipschitz up to some additive error $\varepsilon'_n$ with $\varepsilon'_n \to 0$. 
\item Let $B \subset X$ be an arbitrary Borel subset and let $\mu_{B, n}$ be the probability measure on $X_n$ defined by 
\begin{equation*}
\mu_{B, n} := m_{X_n}(\,\cdot \, \cap p_n^{-1}(B)) / m_{X_n}(p_n^{-1}(B)).
\end{equation*}
Then, for any $\kappa > 0$, we have
\begin{equation*}
\limsup_{n \to \infty}\OD((p_n^{-1}(B), d_{X_n}, \mu_{B, n}); -\kappa) \leq \diam{B}.
\end{equation*}
\item For any two Borel subsets $B_1, B_2 \subset X$ and any $\kappa > 0$, we have
\begin{equation*}
\limsup_{n \to \infty}d_+(p_n^{-1}(B_1), p_n^{-1}(B_2); +\kappa) \leq d_X(B_1, B_2) + \sum_{i=1}^2 \diam{B_i}.
\end{equation*}
\end{enumerate}
\end{thm}

\begin{prop}\label{set_met_prod}
Let $X$ and $Y$ are two metric spaces and let $F$ be a function in $\F^2$. If $F$ is an isotone, then
\begin{equation}
d_F(A \times B, A' \times B') = F(d_X(A, A'), d_Y(B, B'))
\end{equation}
for any subsets $A, A' \subset X$ and $B, B' \subset Y$.
\end{prop}

\begin{proof}
Since $F$ is an isotone, for any $(a, b) \in A \times B$ and $(a', b') \in A' \times B'$,
\begin{equation*}
d_F((a, b), (a', b')) = F(d_X(a, a'), d_Y(b, b')) \geq F(d_X(A, A'), d_Y(B, B')),
\end{equation*}
which implies $d_F(A \times B, A' \times B') \geq F(d_X(A, A'), d_Y(B, B'))$. Let us prove the opposite inequality. We take any two real numbers $\varepsilon, \delta$ such that $d_X(A, A') < \varepsilon$ and $d_Y(B, B') < \delta$. There exist $a \in A$, $a' \in A'$, $b \in B$, and $b' \in B'$ such that $d_X(a, a') < \varepsilon$ and $d_Y(b, b') < \delta$. Then
\begin{equation*}
d_F((a, b), (a', b')) = F(d_X(a, a'), d_Y(b, b')) \leq F(\varepsilon, \delta),
\end{equation*}
which implies $d_F(A \times B, A' \times B') \leq F(\varepsilon, \delta)$. By the continuity of $F$, we have $d_F(A \times B, A' \times B') \leq F(d_X(A, A'), d_Y(B, B'))$. This completes the proof.
\end{proof}

\begin{proof}[Proof of {\rm ``$(\ref{main2}) \Rightarrow (\ref{main1})$''} of Theorem \ref{main}]
The idea of the proof is based on the same of the proof of the fibration theorem.

Assume that the functions $F_n$ satisfy the condition (\ref{main2}) and two sequences $\{X_n\}_{n \in \N}$ and $\{Y_n\}_{n \in \N}$ concentrate to $X$ and $Y$ respectively. Note that the function $F$ is an isotone. By Theorem \ref{equiconc},  there exist Borel measurable maps $p_n \colon X_n \to X$, where $n = 1, 2, \ldots$, that enforce $\varepsilon_n$-concentration of $X_n$ to $X$ and $\prok({p_n}_* m_{X_n}, m_X) \leq \varepsilon_n$ for some sequence $\varepsilon_n \to 0$. Similarly,  there exist Borel measurable maps $q_n \colon Y_n \to Y$, where $n = 1, 2, \ldots$, that enforce $\varepsilon_n$-concentration of $Y_n$ to $Y$ and $\prok({q_n}_*m_{Y_n}, m_Y) \leq \varepsilon_n$. Since
\begin{equation}\label{prodweakconv}
\prok((p_n \times q_n)_*(m_{X_n} \otimes m_{Y_n}), m_X \otimes m_Y) \leq \max\{2\varepsilon_n, 4F(\varepsilon_n, \varepsilon_n)\}
\end{equation}
follows from Lemma \ref{lprok}, it suffices to prove that the map $p_n \times q_n$ enforces $\varepsilon'_n$-concentration of $X_n \times_{F_n} Y_n$ to $X \times_F Y$ for some sequence $\varepsilon'_n \to 0$. By Lemma \ref{mmg5.27} (1) and Lemma \ref{lip_up}, the map $p_n \times q_n$ is $1$-Lipschitz up to $\varepsilon'_n$ for some $\varepsilon'_n \to 0$. By Lemma \ref{mmg5.27} (2), we have
\begin{equation}
(p_n \times q_n)^* \Lip_1(X \times_F Y) \subset U_{2\varepsilon'_n}\left( \Lip_1(X_n \times_{F_n} Y_n)\right).
\end{equation}
Therefore, for any real number $\varepsilon > 0$, it suffices to prove that
\begin{equation}
\Lip_1(X_n \times_{F_n} Y_n) \subset U_{42F(\varepsilon, \varepsilon) + 5\varepsilon}\left( (p_n \times q_n)^* \Lip_1(X \times_F Y)\right)
\end{equation}
holds for every sufficiently large $n$.

We take any $\varepsilon > 0$ and any $f_n \in \Lip_1(X_n \times_{F_n} Y_n)$. There are finitely many mutually disjoint nonempty open subsets $B_1^X, B_2^X, \ldots, B_N^X \subset X$ such that $m_X(\partial B_i^X) = 0$, $\diam{B_i^X} < \varepsilon$, and
\begin{equation*}
m_X\left( X \setminus \bigcup_{i = 1}^N B_i^X \right) < \varepsilon.
\end{equation*}
Similarly, there exist mutually disjoint nonempty open subsets $B_1^Y, B_2^Y,$ $\ldots, B_N^Y \subset Y$ such that $m_Y(\partial B_i^Y) = 0$, $\diam{B_i^Y} < \varepsilon$, and
\begin{equation*}
m_Y\left( Y \setminus \bigcup_{i = 1}^N B_i^Y \right) < \varepsilon.
\end{equation*} For each $i = 1, 2, \ldots, N$, we take points $x_i \in B_i^X$ and $y_i \in B_i^Y$ and fix them. We put
\begin{align*}
A_{in}^X := p_n^{-1}(B_i^X), \quad  & \mu_{in} := {m_{X_n}(A_{in}^X)}^{-1} m_{X_n}|_{A_{in}^X}, \\
A_{in}^Y := q_n^{-1}(B_i^Y), \quad &  \nu_{in} := {m_{Y_n}(A_{in}^Y)}^{-1} m_{Y_n}|_{A_{in}^Y}
\end{align*}
for every $i$. Note that $m_{X_n}(A_{in}^X)$ converges to $m_X(B_i^X)$ and $m_{Y_n}(A_{in}^Y)$ converges to $m_Y(B_i^Y)$ as $n \to \infty$. We define a Borel measurable map $g_n \colon X \times Y \to \R$ by
\begin{equation}
g_n(x, y) := \left\{ \begin{array}{ll} \displaystyle \lm(f_n |_{A_{in}^X \times A_{jn}^Y} ; \mu_{in} \otimes \nu_{jn}) &  \text{if } (x, y) \in B_i^X \times B_j^Y, \\ \displaystyle  0 & \text{otherwise.} \end{array} \right.
\end{equation}
Our immediate goal is to prove that $g_n$ is $1$-Lipschitz up to $30F(\varepsilon, \varepsilon) + 2\varepsilon$ with respect to $d_F$ and $(p_n \times q_n)_* (m_{X_n} \otimes m_{Y_n})$ for every sufficiently large $n$.

Setting $\rho_{ijkl} := F(d_X(x_i, x_j), d_Y(y_k, y_l)) + 4F(\varepsilon, \varepsilon) + \varepsilon$ for any $i, j, k, l = 1, 2, \ldots, N$, we find $\lambda > 0$ such that $0 < \lambda \rho_{ijkl} < 1 / 4$. 

\begin{claim}\label{lprok_claim}
For each $i, j, k, l$ and every sufficiently large $n$, we have
\begin{equation}\label{lprok_claim_eq}
\lprok{\lambda}(\mu_{in} \otimes \nu_{kn}, \mu_{jn} \otimes \nu_{ln}) \leq \rho_{ijkl},
\end{equation}
where $\lprok{\lambda}$ in the left-hand side is with respect to $d_{F_n}$.
\end{claim}

\begin{proof}
The proof of the claim is similar to that of \cite{MMG}*{Claim 9.9}.

We fix $i, j, k$, and $l$. By the $\pi$-$\lambda$ theorem, it is sufficient to prove that
\begin{equation}\label{lprok_claim_goal}
\mu_{jn} \otimes \nu_{ln}(U_{\rho_{ijkl}}(C_n \times D_n)) \geq \mu_{in} \otimes \nu_{kn}(C_n \times D_n) - \lambda \rho_{ijkl}
\end{equation}
for any Borel subset $C_n \subset X_n$ and $D_n \subset Y_n$. Take any Borel subsets $C_n \subset X_n$ and $D_n \subset Y_n$. We are able to assume that $C_n \subset A_{in}^X$ and $D_n \subset A_{kn}^Y$ since $\mu_{in}(C_n) = \mu_{in}(C_n \cap A_{in}^X)$ and $\nu_{kn}(D_n) = \nu_{kn}(D_n \cap A_{kn}^Y)$. Let $\kappa$ be a real number such that
\begin{equation*}
0 < \kappa \leq \lambda \rho_{ijkl} \inf_{n \in \N} \min\{m_{X_n}(A_{in}^X), m_{X_n}(A_{jn}^X), m_{Y_n}(A_{kn}^Y), m_{Y_n}(A_{ln}^Y)\}. 
\end{equation*}
If $m_{X_n}(C_n) < \kappa$ or $m_{Y_n}(D_n) < \kappa$, then we have
\begin{equation*}
\mu_{in} \otimes \nu_{kn}(C_n \times D_n) = \frac{m_{X_n}(C_n) \, m_{Y_n}(D_n)}{m_{X_n}(A_{in}^X) \, m_{Y_n}(A_{kn}^Y)} \leq \lambda \rho_{ijkl}, 
\end{equation*}
so that we obtain \eqref{lprok_claim_goal}. Assume that $m_{X_n}(C_n) \geq \kappa$ and $m_{Y_n}(D_n) \geq \kappa$. We define two functions $\varphi_n \colon A_{jn}^X \to \R$ and $\psi_n \colon A_{kn}^Y \to \R$ by
\begin{equation*}
\varphi_n(x) := d_{X_n}(x, C_n) \text{ and } \psi_n(y) := d_{Y_n}(y, D_n)
\end{equation*}
for $x \in A_{jn}^X$ and $y \in A_{ln}^Y$, and let
\begin{align*}
E_n^X & := \left\{x \in A_{jn}^X \midd |\varphi_n(x) - \lm(\varphi_n ; \mu_{jn})| \leq \varepsilon \right\}, \\
E_n^Y & := \left\{y \in A_{ln}^Y \midd |\psi_n(x) - \lm(\psi_n ; \nu_{ln})| \leq \varepsilon \right\}.
\end{align*}
For any $\kappa' \in (0, 1/2)$ and every sufficiently large $n$, by Theorem \ref{Fib}, we have
\begin{align*}
\OD(\mu_{jn}; -\kappa') & \left( := \OD((A_{jn}^X, \mu_{jn}); -\kappa') \right) < \varepsilon, \\ 
\OD(\nu_{ln}; -\kappa') & \left( := \OD((A_{ln}^Y, \mu_{ln}); -\kappa') \right) < \varepsilon,
\end{align*}
and then, by Lemma \ref{mmg7.31}, we have
\begin{equation*}
\LR(\mu_{jn}; -\kappa') < \varepsilon \text{ and } \LR(\nu_{ln}; -\kappa') < \varepsilon.
\end{equation*}
Thus we have $\mu_{jn}(E_n^X)$, $\nu_{ln}(E_n^Y) \to 1$ as $n \to \infty$, which imply
\begin{equation*}
m_{X_n}(E_n^X) \geq \kappa, \ m_{Y_n}(E_n^Y) \geq \kappa, \text{ and } \mu_{jn}\otimes \nu_{ln}(E_n^X \times E_n^Y) \geq 1 - \lambda \rho_{ijkl}
\end{equation*}
for every sufficiently large $n$. By Theorem \ref{Fib}, it holds that 
\begin{equation*}
d_{X_n}(C_n, E_n^X) \leq d_+(A_{in}^X, A_{jn}^X; +\kappa) < d_X(B_i, B_j) + 2\varepsilon \leq d_X(x_i, x_j) + 2\varepsilon
\end{equation*}
for every sufficiently large $n$. For any two points $x, x' \in E_n^X$, we have
\begin{equation*}
d_{X_n}(x, C_n) \leq \lm(\varphi_n; \mu_{jn}) + \varepsilon \leq d_{X_n}(x', C_n) + 2\varepsilon,
\end{equation*}
which implies $d_{X_n}(x, C_n) \leq d_{X_n}(E_n^X, C_n) + 2\varepsilon < d_X(x_i, x_j) + 4\varepsilon$ for every $x \in E_n^X$. Similarly, $d_{Y_n}(y, D_n) < d_Y(y_k, y_l) + 4\varepsilon$ also holds for every $y \in E_n^Y$. By Proposition \ref{set_met_prod}, we have
\begin{align*}
& d_{F_n}((x, y), C_n \times D_n) \\
& = F_n(d_{X_n}(x, C_n), d_{Y_n}(y, D_n)) \\
& < F(d_{X_n}(x, C_n), d_{Y_n}(y, D_n)) + \varepsilon \\
& \leq F(d_X(x_i, x_j) + 4\varepsilon, d_Y(y_k, y_l) + 4\varepsilon) + \varepsilon \\
& \leq F(d_X(x_i, x_j), d_Y(y_k, y_l)) + 4F(\varepsilon, \varepsilon) + \varepsilon = \rho_{ijkl}
\end{align*}
for any $(x, y) \in E_n^X \times E_n^Y$ and every sufficiently large $n$, where the first inequality follows from that $F_n$ uniformly converges to $F$ on any compact sets, and the second follows from that $F$ is an isotone. This means $E_n^X \times E_n^Y \subset U_{\rho_{ijkl}}(C_n \times D_n)$. Therefore we have
\begin{equation*}
\begin{split}
& \mu_{jn} \otimes \nu_{ln}\left(U_{\rho_{ijkl}}(C_n \times D_n)\right) \geq \mu_{jn} \otimes \nu_{ln}(E_n^X \times E_n^Y) \\
& \geq 1 - \lambda \rho_{ijkl} \geq  \mu_{in} \otimes \nu_{kn}(C_n \times D_n) - \lambda \rho_{ijkl},
\end{split}
\end{equation*}
so that we obtain \eqref{lprok_claim_goal}. This completes the proof.
\end{proof}

By Claim \ref{lprok_claim} and Strassen's theorem (Theorem \ref{Strassen}), there exists $\rho_{ijkl}$-subtransport plan $\pi_{ijkl}^n$ between $\mu_{in} \otimes \nu_{kn}$ and $\mu_{jn} \otimes \nu_{ln}$ such that $\df{\pi_{ijkl}^n} \leq \lambda \rho_{ijkl}$. Since $\df{\pi_{ijkl}^n} < 1/4$, we have
\begin{align*}
&|g_n(x, y) - g_n(x', y')| \\
& = |\lm(f_n |_{A_{in}^X \times A_{kn}^Y} ; \mu_{in} \otimes \nu_{kn} ) - \lm(f_n |_{A_{jn}^X \times A_{ln}^Y} ; \mu_{jn} \otimes \nu_{ln} )| \\
& \leq \rho_{ijkl} + \OD(\mu_{in} \otimes \nu_{kn} ; -\frac{1}{4}) + \OD(\mu_{jn} \otimes \nu_{ln} ; -\frac{1}{4})
\end{align*}
for any $(x, y) \in B_i^X \times B_k^Y$ and any $(x', y') \in B_j^X \times B_l^Y$, where the last inequality follows from Lemma \ref{lm_lem}. We have
\begin{align*}
& \rho_{ijkl} = F(d_X(x_i, x_j), d_Y(y_k, y_l)) + 4F(\varepsilon, \varepsilon) + \varepsilon \\
& \leq F(d_X(x, x') + 2\varepsilon, d_Y(y, y') + 2\varepsilon)+ 4F(\varepsilon, \varepsilon) + \varepsilon \\
& \leq d_F((x, y), (x', y')) + 6F(\varepsilon, \varepsilon) + \varepsilon,
\end{align*}
where the first inequality follows from that $F$ is an isotone. Moreover, by Lemma \ref{key_F}, we have
\begin{equation*}
\begin{split}
& \OD(\mu_{in} \otimes \nu_{kn} ; -\frac{1}{4}) \\
& \leq 4F(\OD(\mu_{in} ; - \frac{1}{16}), 0) + 8F(0, \OD(\nu_{kn} ; -\frac{1}{16})) \\
& \leq 12F(\OD(\mu_{in} ; - \frac{1}{16}), \OD(\nu_{kn} ; -\frac{1}{16})),
\end{split}
\end{equation*}
and, by Theorem \ref{Fib}, we see that
\begin{equation*}
\OD(\mu_{in} ; - \frac{1}{16}) < \varepsilon \text{ and } \OD(\nu_{kn} ; -\frac{1}{16}) < \varepsilon
\end{equation*}
for every sufficiently large $n$. Thus we obtain
\begin{equation}\label{lip_up_eq}
|g_n(x, y) - g_n(x', y')| \leq d_F((x, y), (x', y')) + 30F(\varepsilon, \varepsilon) + \varepsilon
\end{equation}
for any $(x, y), (x', y') \in \bigcup_{i, j = 1}^N B_i^X \times B_j^Y$ and every sufficiently large $n$. Furthermore it holds that
\begin{align*}
&\lim_{n \to \infty} (p_n \times q_n)_* (m_{X_n} \otimes m_{Y_n}) \left( \bigcup_{i,j=1}^N B_i^X \times B_j^Y \right) \\
& = m_X \otimes m_Y \left( \bigcup_{i,j=1}^N B_i^X \times B_j^Y \right) = m_X \left( \bigcup_{i=1}^N B_i^X \right) m_Y \left( \bigcup_{j=1}^N B_j^Y \right) \\
& \geq (1 - \varepsilon)^2 > 1 - 2\varepsilon,
\end{align*}
where the first equality follows from (\ref{prodweakconv}). Combining this and (\ref{lip_up_eq}) implies that $g_n$ is $1$-Lipschitz up to $30F(\varepsilon, \varepsilon) + 2\varepsilon$ with respect to $d_F$ and $(p_n \times q_n)_* (m_{X_n} \otimes m_{Y_n})$ for every sufficiently large $n$. By Lemma \ref{MHext}, we see that there exists $\tilde{g}_n \in \Lip_1(X \times Y)$ such that
\begin{equation*}
\kf^{m_{X_n} \otimes m_{Y_n}}\left((p_n \times q_n)^* g_n, (p_n \times q_n)^* \tilde{g}_n \right) \leq 30F(\varepsilon, \varepsilon) + 2\varepsilon.
\end{equation*}

Let $\kappa := \min{\{ \varepsilon/ N^2 , 1/4\}}$. For every sufficiently large $n$, we have
\begin{equation*}
\LR(\mu_{in}\otimes \nu_{jn} ; -\kappa) \leq \OD(\mu_{in}\otimes \nu_{jn} ; -\kappa) \leq 12F(\varepsilon, \varepsilon)
\end{equation*}
for each $i, j = 1, 2, \ldots, N$. Setting 
\begin{equation*}
K_{ijn} := \left\{(x, y) \midd |f_n(x, y) - \lm(f_n |_{A_{in}^X \times A_{jn}^Y} ; \mu_{in} \otimes \nu_{jn} )| > 12F(\varepsilon, \varepsilon)\right\},
\end{equation*}
we have
\begin{align*}
& m_{X_n} \otimes m_{Y_n} (|f_n - (p_n \times q_n)^* g_n| > 12F(\varepsilon, \varepsilon)) \\
 \leq & \, \sum_{i, j = 1}^N m_{X_n} \otimes m_{Y_n}\left((A_{in}^X  \times A_{jn}^Y) \cap K_{ijn} \right) \\
& + m_{X_n} \otimes m_{Y_n} \left((X_n \times Y_n) \setminus  \bigcup_{i,j=1}^N A_{in}^X \times A_{jn}^Y \right) \\
 \leq & \, N^2 \kappa + (p_n \times q_n)_* (m_{X_n} \otimes m_{Y_n}) \left( (X \times Y) \setminus \bigcup_{i,j=1}^N B_i^X \times B_j^Y \right) \\
 \leq & \, \varepsilon + 2\varepsilon = 3\varepsilon.
\end{align*}
Thus we obtain
\begin{equation*}
\kf^{m_{X_n} \otimes m_{Y_n}}\left(f_n, (p_n \times q_n)^* g_n \right) < 12F(\varepsilon, \varepsilon) + 3\varepsilon,
\end{equation*}
which implies
\begin{equation*}
\kf^{m_{X_n} \otimes m_{Y_n}}\left(f_n, (p_n \times q_n)^* \tilde{g}_n \right) < 42F(\varepsilon, \varepsilon) + 5\varepsilon.
\end{equation*}
This completes the proof.
\end{proof}

\begin{proof}[Proof of Corollary \ref{main_cor}]
We just apply Theorem \ref{main} with $F_n = F_{p_n}$, $n = 1, 2,\ldots$, and $F = F_p$, where $F_p$ is the function of Example \ref{exF} (\ref{exF1}). 
\end{proof}

\begin{proof}[Proof of $(\ref{main_1dim2}) \Rightarrow (\ref{main_1dim1})$ of Theorem \ref{main_1dim}]
Let $Y_n$, $n = 1, 2, \ldots$, and $Y$ be one-point mm-spaces and let $G_n$ and $G$ be the functions on $[0, +\infty)^2$ defined by
\begin{equation*}
G_n(s, t) := F_n(s) + t \text{ and } G(s, t) := F(s) + t
\end{equation*}
for $s, t \in [0, + \infty)$. We apply the implication from (\ref{main2}) to (\ref{main1}) of Theorem \ref{main}. 
\end{proof}

\subsection{A new specific example of the concentration}

\begin{ex}\label{5.1}
We consider the $n$-dimensional unit sphere $S^n(1)$ and  the interval $[0, \pi]$. These spaces are both equipped with the distance and normalized measure induced by the standard Riemannian metric. We take an arbitrary point $\bar{x} \in S^n(1)$ and fix it. We attach the interval $[0, \pi]$ to the sphere $S^n(1)$ at their points $\pi \in [0, \pi]$ and $\bar{x} \in S^n(1)$, and denote their united space by $X_n$. That is, the space $X_n$ is defined as the mm-space
\begin{equation*}
X_n := \fracinline{ [0, \pi] \sqcup S^n(1) }{\pi = \bar{x}},
\end{equation*}
where the distance $d_{X_n}$ is defined by
\begin{equation*}
d_{X_n}(x, x') := \left\{\begin{array}{ll} d_{[0, \pi]}(x, x') & \text{ if } x, x' \in [0, \pi], \\ d_{S^n(1)}(x, x') & \text{ if } x, x' \in S^n(1), \\  d_{[0, \pi]}(x, \pi) + d_{S^n(1)}(\bar{x}, x') & \text{ if } x \in [0, \pi], x' \in S^n(1), \\ d_{[0, \pi]}(x', \pi) + d_{S^n(1)}(\bar{x}, x) & \text{ if } x' \in [0, \pi], x \in S^n(1) \end{array}\right.
\end{equation*}
for $x, x' \in X_n$, and the measure $m_{X_n}$ is defined by
\begin{equation*}
m_{X_n} := \frac{1}{2} m_{[0, \pi]} + \frac{1}{2} m_{S^n(1)}.
\end{equation*}
The sequence $\{X_n\}_{n \in \N}$ concentrates to the following mm-space $X$. The mm-space $X$ is the subset $[0, \pi] \cup \{3\pi/2\}$ of the one-dimensional Euclidean space
$\R$ with the Euclidean distance and the measure
\begin{equation*}
m_{X} := \frac{1}{2} m_{[0, \pi]} + \frac{1}{2} \delta_{\frac{3}{2}\pi},
\end{equation*}
where $\delta_x$ is the Dirac measure at a point $x$. This is proved by applying Theorem \ref{equiconc} to the maps $p_n \colon X_n \to X$, $n = 1, 2, \ldots$, defined by
\begin{equation}
p_n(x) := \left\{\begin{array}{ll} x & \text{ if } x \in [0, \pi], \\ \frac{3}{2}\pi & \text{ if } x \in S^n(1). \end{array}\right.
\end{equation}
If the reader wishes to prove its details, one reads Section 5. The proof is similar to that of Claim \ref{claim_nece_1dim} and \ref{claim_nece} (and is easier than them).
\end{ex}

Applying the implication from (\ref{main2}) to (\ref{main1}) of Theorem \ref{main}, we understand the concentration of product spaces of two copies of $X_n$ in Example \ref{5.1}.

\begin{ex}\label{5.2}
Let $X_n$, $n = 1, 2, \ldots$, and $X$ be mm-spaces of Example \ref{5.1}. Corollary \ref{main_cor} implies that the sequence of the $l_p$-product spaces $\{X_n  \times_p X_n\}_{n \in \N}$ concentrates to the $l_p$-product space $X \times_p X$ for any $p \in [1, +\infty]$. The limit space $X \times_p X$ is mm-isomorphic to the subset
\begin{equation*}
\left\{ (x, y) \in \R^2 \mid x, y \in [0, \pi] \cup \{3\pi/2\} \right\}
\end{equation*}
of the $l_p$-normed space $(\R^2, \|\cdot \|_p)$.
\end{ex}

\section{The necessity of the isotonicity}
In this section, we prove the implication from (\ref{main1}) to (\ref{main2}) of Theorem \ref{main} and Theorem \ref{main_1dim}. In order to prove them, we construct some counterexample of the condition (\ref{main1}) if the condition (\ref{main2}) does not hold.

We first prove Theorem \ref{main_1dim}.

\begin{proof}[Proof of {\rm ``$(\ref{main_1dim1}) \Rightarrow (\ref{main_1dim2})$''} of Theorem \ref{main_1dim}]
Assume that the condition (\ref{main_1dim2}) does not hold. That is, up to choosing a subsequence of $n$, we are able to assume that there exist two real numbers $s, \eta > 0$ and a sequence $\{s_n\}_{n \in \N} \subset (0, +\infty)$ such that 
\begin{equation*}
s < s_n \text{ and } F_n(s) > F_n(s_n) + \eta
\end{equation*}
for any $n \in \N$. Moreover, we can assume that 
\begin{equation}\label{s'_min}
F_n(s_n) = \min_{s \leq t \leq s_n} F_n(t).
\end{equation}
Choosing a subsequence of $n$, we can assume that there exists limit of $\{F_n(s_n)\}_{n \in \N}$ as $n \to \infty$. we see that
\begin{equation*}
F(s) \geq \lim_{n \to \infty} F_n(s_n) + \eta.
\end{equation*}
We define an mm-space $X$ as
\begin{equation*}
X := (\{x_0, x_1\}, d_X, \frac{1}{2} \delta_{x_0} + \frac{1}{2} \delta_{x_1}), \quad d_X(x_0, x_1) := s.
\end{equation*}
We set
\begin{equation*}
r_n := \frac{\sqrt{s_n^2 - s^2}}{2} > 0 \text{ and } k_n := \max\{n, \lceil r_n^4\rceil\}
\end{equation*}
for each $n$, where $\lceil\cdot\rceil$ is the ceiling function. Let $S^{k_n}(r_n)$ be the $k_n$-dimensional sphere of radius $r_n$ in $\R^{k_n + 1}$ centered at the origin. The sphere $S^{k_n}(r_n)$ is equipped with the Euclidean distance $\|\cdot \|$ and the normalized probability volume measure $\sigma^{k_n}$. Define an mm-space $X_n$ for each $n$ as 
\begin{equation*}
X_n := X \times_2 (S^{k_n}(r_n), \|\cdot \|, \sigma^{k_n}).
\end{equation*}
Note that embedding $X$ into the 1-dimensional Euclidean space, $X_n$ is regarded as a subset of the $(k_n+2)$-dimensional Euclidean space $(\R^{k_n+2}, \|\cdot\|)$ naturally. By Corollary \ref{main_cor} and Example \ref{sphere_norm}, the sequence $\{X_n\}_{n \in \N}$ concentrates to $X$. Let us prove the following claim.

\begin{claim}\label{claim_nece_1dim}
The sequence $\{(X_n, F_n \circ d_{X_n}, m_{X_n})\}_{n \in \N}$ concentrates to the mm-space $Y$ defined by
\begin{equation*}
Y := (\{y_0, y_1\}, d_Y, \frac{1}{2} \delta_{y_0} + \frac{1}{2} \delta_{y_1}), \quad d_Y(y_0, y_1) := \lim_{n \to \infty} F_n(s_n),
\end{equation*} 
as $n \to \infty$.
\end{claim}

\begin{proof}
For $i = 0, 1$, we set a subset $S^n_i$ of $X_n$ and a measure $\sigma^n_i$ on $S^n_i$ by
\begin{equation*}
S^n_i := \{x_i\} \times S^{k_n}(r_n) \subset X_n \subset \R^{k_n+2}, \quad \sigma^n_i := \delta_{x_i} \otimes \sigma^{k_n}.
\end{equation*}
Note that $(S^n_i, \|\cdot\|, \sigma^n_i)$ is mm-isomorphic to $(S^{k_n}(r_n), \|\cdot\|, \sigma^{k_n})$ for both $i = 0, 1$. Let $p_n \colon X_n \to Y$ be the map defined by
\begin{equation*}
p_n(x) := \left\{ \begin{array}{ll} y_0 & \text{ if } x \in S^n_0, \\ y_1 & \text{ if } x \in S^n_1. \end{array} \right.
\end{equation*}
Note that ${p_n}_* m_{X_n} = m_Y$. Let $\varepsilon > 0$ be a sufficiently small arbitrary real number. We find a number $N \in \N$ such that 
\begin{align*}
&|F_n(\varepsilon) - F(\varepsilon)| < \varepsilon, \quad |F_n(s_n) - d_Y(y_0, y_1)| < \varepsilon, \\
& \text{and } \OD((S^n_i, \| \cdot \|, \sigma^n_i); -\varepsilon) < \varepsilon 
\end{align*}
hold for any $n \geq N$ and $i = 0, 1$. Let us prove that $p_n$ enforces $(24F(\varepsilon) + 27\varepsilon)$-concentration of $(X_n, F_n \circ d_{X_n}, m_{X_n})$ to $Y$ for any $n \geq N$. We fix $n$ with $n \geq N$. For any $x \in S^n_0$ and $x' \in S^n_1$, it holds that
\begin{equation*}
s \leq \|x - x'\| \leq s_n,
\end{equation*}
so that
\begin{equation*}
F_n(\|x - x'\|) \geq F_n(s_n) > d_Y(y_0, y_1) - \varepsilon,
\end{equation*}
where the first inequality follows from \eqref{s'_min}. Thus, the map $p_n$ is 1-Lipschitz up to $\varepsilon$ with respect to $F_n \circ d_{X_n}$. By Lemma \ref{mmg5.27}, we have
\begin{equation*}
{p_n}^* \Lip_1(Y) \subset U_{2\varepsilon}(\Lip_1(X_n)).
\end{equation*}
We prove the other side inclusion. We take any function $f_n \in \Lip_1(X_n, {F_n \circ d_{X_n}})$ and define a function $g_n : Y \to \R$ by
\begin{equation*}
g_n(y_0) := \lm(f_n; \sigma^n_0), \quad  g_n(y_1) := \lm(f_n; \sigma^n_1).
\end{equation*}
By Lemma \ref{mmg7.31} and Lemma \ref{key_1dim},
\begin{align*}
& \LR((S^n_i, F_n \circ \| \cdot \|, \sigma^n_i); -2\varepsilon) \\
& \leq \OD((S^n_i, F_n \circ \| \cdot \|, \sigma^n_i); -2\varepsilon) \\
& \leq 4F_n(\OD((S^n_i, \| \cdot \|, \sigma^n_i); -\varepsilon)) \leq 8F_n(\varepsilon).
\end{align*}
Thus we have 
\begin{align*}
& m_{X_n}\left\{x \in X_n \midd |f_n(x) - {p_n}^* g_n(x)| > 8F_n(\varepsilon) \right\} \\
& = \frac{1}{2} \sum_{i=0}^1  \sigma^n_i \left\{x \in S^n_i \midd |f_n(x) - \lm(f_n; \sigma^n_i)| > 8F_n(\varepsilon) \right\} \leq 2\varepsilon,
\end{align*}
which implies $\kf^{m_{X_n}}(f_n, {p_n}^*g_n) < 8F_n(\varepsilon) + 2\varepsilon < 8F(\varepsilon) + 10\varepsilon$. Let $T_n \colon S^n_0 \to S^n_1$ be the map defined by
\begin{equation*}
T_n(x_0, a) := (x_1, -a)
\end{equation*}
for $a \in S^{k_n}(r_n)$. Note that ${T_n}_* \sigma^n_0 = \sigma^n_1$. For any $x \in S^n_0$, we have 
\begin{equation*}
\|x - T_n(x)\|^2 = s^2 + (2r_n)^2 = s_n^2,
\end{equation*}
so that
\begin{equation*}
F_n(\|x - T_n(x) \|) = F_n(s_n).
\end{equation*}
Thus, the measure $(\id, T_n)_* \sigma^n_0$ is an $F_n(s_n)$-(sub)transport plan between $\sigma^n_0$ and $\sigma^n_1$ (with $\df{((\id, T_n)_* \sigma^n_0)} = 0$). By Lemma \ref{lm_lem} and Lemma \ref{key_1dim}, we have
\begin{align*}
& |\lm(f_n; \sigma^n_0) - \lm(f_n; \sigma^n_1)|  \\
& \leq F_n(s_n) + \sum_{i = 0}^1 \OD((S^n_i, F_n \circ \|\cdot\|, \sigma^n_i); -2\varepsilon) \\
& \leq F_n(s_n) + \sum_{i = 0}^1 4F_n(\OD((S^n_i, \|\cdot\|, \sigma^n_i); -\varepsilon)) \\
& \leq F_n(s_n)+ 16F_n(\varepsilon) < d_Y(y_0, y_1) + 16F(\varepsilon) + 17\varepsilon.
\end{align*}
Lemma \ref{MHext} implies $\kf^{m_X}(g_n, \Lip_1(Y)) < 16F(\varepsilon) + 17\varepsilon$. Taking
\begin{equation*}
\kf^{m_Y}(g_n, \Lip_1(Y)) =  \kf^{m_{X_n}}({p_n}^*g_n, {p_n}^* \Lip_1(Y))
\end{equation*}
into account, we have
\begin{equation*}
\begin{split}
&\kf^{m_{X_n}}(f_n, {p_n}^*\Lip_1(Y)) \\ 
&\leq \kf^{m_{X_n}}(f_n, {p_n}^* g_n) + \kf^{m_Y}(g_n, \Lip_1(Y)) < 24F(\varepsilon) + 27\varepsilon.
\end{split}
\end{equation*}
Thus we obtain $\Lip_1(X_n) \subset U_{24F(\varepsilon) + 27\varepsilon}({p_n}^* \Lip_1(Y))$ and then the map $p_n$ enforces $(24F(\varepsilon) + 27\varepsilon)$-concentration of $(X_n, F_n \circ d_{X_n}, m_{X_n})$ to $Y$ for every $n \geq N$. By Theorem \ref{equiconc}, the sequence $\{(X_n, F_n \circ d_{X_n}, m_{X_n})\}_{n \in \N}$ concentrates to $Y$ as $n \to \infty$. The proof of the claim is now completed.
\end{proof}
Since
\begin{equation*}
d_X(x_0, x_1) = F(s) \geq \lim_{n \to \infty} F_n(s_n) + \eta = d_Y(y_0, y_1) + \eta,
\end{equation*}
the mm-space $(X, F \circ d_X, m_X)$ is not mm-isomorphic to $Y$. Thus Claim \ref{claim_nece_1dim} means that the condition (\ref{main_1dim1}) of Theorem \ref{main_1dim} does not hold. Therefore we obtain the implication from (\ref{main_1dim1}) to (\ref{main_1dim2}) of Theorem \ref{main_1dim}.
\end{proof}

We next prove Theorem \ref{main}. The idea of the proof is same as that of the above proof but the following proof is more complicated.

\begin{proof}[Proof of {\rm ``$(\ref{main1}) \Rightarrow (\ref{main2})$''} of Theorem \ref{main}]
Assume that the condition (\ref{main2}) does not hold. Up to choosing a subsequence of $n$, we are able to assume that there exist a real number $\eta > 0$, a pair $(s, t) \in [0, +\infty)^2$, and a sequence $\{(s_n, t_n)\}_{n \in \N} \subset [0, +\infty)^2$ such that 
\begin{equation*}
s < s_n, \ t < t_n, \text{ and } F_n(s, t) > F_n(s_n, t_n) + \eta
\end{equation*}
for any $n \in \N$. We define two mm-spaces $X$ and $Y$ as
\begin{align*}
X := & \ (\{\bar{x}_0, \bar{x}_1\}, d_X, \frac{1}{2} \delta_{\bar{x}_0} + \frac{1}{2} \delta_{\bar{x}_1}), & d_X(\bar{x}_0, \bar{x}_1) := s, \\
Y := & \ (\{\bar{y}_0, \bar{y}_1\}, d_Y, \frac{1}{2} \delta_{\bar{y}_0} + \frac{1}{2} \delta_{\bar{y}_1}), & d_Y(\bar{y}_0, \bar{y}_1) := t.
\end{align*}
For each $n$, let
\begin{align*}
r_n := \frac{\sqrt{s_n^2 - s^2}}{2} > 0, \quad & \rho_n := \frac{\sqrt{t_n^2 - t^2}}{2} > 0, \\
k_n := 2\max\{n, \lceil r_n^4\rceil\} + 1, \quad & l_n := 2\max\{n, \lceil \rho_n^4\rceil\} + 1.
\end{align*}
Define two mm-spaces $X_n$ and $Y_n$ for each $n$ as 
\begin{equation*}
\begin{split}
X_n &  := X \times_2 (S^{k_n}(r_n), \|\cdot \|, \sigma^{k_n}), \\
Y_n &  := Y \times_2 (S^{l_n}(\rho_n), \|\cdot \|, \sigma^{l_n}).
\end{split}
\end{equation*}
Note that we regard $X_n$ and $Y_n$ as subsets of the Euclidean spaces with dimensions $k_n+2$ and $l_n+2$ respectively, and that the two sequences $\{X_n\}_{n \in \N}$ and $\{Y_n\}_{n \in \N}$ concentrate to $X$ and $Y$ respectively.

In addition, we define three numbers $\alpha_n$, $\beta_n$, $\gamma_n$ for each $n$ by
\begin{align*}
\alpha_n & := \min_{\substack{s \leq u_1 \leq s_n, \\ 0 \leq v_1 \leq 2\rho_n}} F_n(u_1, v_1), \\ 
\beta_n & := \min_{\substack{0 \leq u_2 \leq 2r_n, \\ t \leq v_2 \leq t_n}} F_n(u_2, v_2), \\ 
\gamma_n & := \min_{\substack{s \leq u_3 \leq s_n \\ t \leq v_3 \leq t_n}} F_n(u_3, v_3).
\end{align*}

\begin{claim}\label{claim_Z}
For each $n \in \N$, the triplet $(\alpha_n, \beta_n, \gamma_n)$ is a triangle triplet.
\end{claim}

\begin{proof}
We fix $n \in \N$ and take any $u_1, u_3 \in [s, s_n]$, $u_2 \in [0, 2r_n]$, $v_2, v_3 \in [t, t_n]$, $v_1 \in [0, 2\rho_n]$. We first prove $\alpha_n \leq \beta_n + \gamma_n$. If $u_2 \leq u_3$, then we have
\begin{equation*}
\alpha_n \leq F_n(u_3, |v_2 - v_3|) \leq F_n(u_2, v_2) +  F_n(u_3, v_3)
\end{equation*}
since $|v_2 - v_3| \leq t_n - t \leq 2\rho_n$ and $(u_3, u_2, u_3)$, $(|v_2 - v_3|, v_2, v_3)$ are triangle triplets. If $u_2 \geq u_3$, then we see that
\begin{equation*}
s \leq u_3 \leq u_2 \leq 2r_n \leq s_n.
\end{equation*}
By this, we have
\begin{equation*}
\alpha_n \leq F_n(u_2, |v_2 - v_3|) \leq F_n(u_2, v_2) +  F_n(u_3, v_3).
\end{equation*}
These inequalities imply $\alpha_n \leq \beta_n + \gamma_n$. We have $\beta_n \leq \alpha_n + \gamma_n$ by the symmetric discussion as the proof of $\alpha_n \leq \beta_n + \gamma_n$. We next prove $\gamma_n \leq \alpha_n + \beta_n$. If $u_{3-i} \leq u_i$ and $v_{3-j} \leq v_j$ for $i, j = 1, 2$, then we have 
\begin{equation*}
s \leq u_i \leq s_n, \quad t \leq v_j \leq t_n
\end{equation*}
and then
\begin{equation*}
\gamma_n \leq F_n(u_i, v_j) \leq F_n(u_1, v_1) +  F_n(u_2, v_2),
\end{equation*}
which implies $\gamma_n \leq \alpha_n + \beta_n$. The proof is completed.
\end{proof}

Choosing a subsequence of $n$, we can assume that there exist limits of $\{\alpha_n\}_{n \in \N}$, $\{\beta_n\}_{n \in \N}$, and $\{\gamma_n\}_{n \in \N}$ as $n \to \infty$ and we denote these limits by $\alpha$, $\beta$, and $\gamma$ respectively. Note that $\alpha$, $\beta$, and $\gamma$ are positive, their triplet $(\alpha, \beta, \gamma)$ is a triangle triplet, and
\begin{equation*}
\alpha \leq F(s, 0), \quad \beta \leq F(0, t), \quad \gamma \leq F(s, t) - \eta.
\end{equation*}
In fact, for each $n$, it holds that 
\begin{equation*}
\alpha_n \leq F_n(s, 0), \quad \beta_n \leq F_n(0, t), \quad \gamma_n \leq F_n(s_n, t_n) \leq F_n(s, t) - \eta.
\end{equation*}
Define an mm-space $Z$ by
\begin{equation*}
Z := (\{z_{00}, z_{10}, z_{01}, z_{11}\}, d_Z, \frac{1}{4}\sum^1_{i,j = 0} \delta_{z_{ij}}),
\end{equation*}
where $d_Z$ is a metric on $Z$ defined as
\begin{equation*}
\begin{split}
d_Z(z_{ij}, z_{ij}) := 0, & \quad d_Z(z_{ij}, z_{1-i, j}) := \alpha, \\ d_Z(z_{ij}, z_{i,1-j})  := \beta, & \quad d_Z(z_{ij}, z_{1-i, 1-j})  := \gamma
\end{split}
\end{equation*}
for every $i, j = 0, 1$. Let us prove the following claim.

\begin{claim}\label{claim_nece}
The sequence $\{X_n \times_{F_n} Y_n\}_{n \in \N}$ concentrates to $Z$ as $n \to \infty$.
\end{claim}

\begin{proof}
Let
\begin{align*}
S^n_{i} := \{\bar{x}_i\} \times S^{k_n}(r_n) \subset X_n, & \quad \sigma^n_i := \delta_{\bar{x}_i} \otimes \sigma^{k_n}, \\
T^n_{j} := \{\bar{y}_j\} \times S^{l_n}(\rho_n) \subset Y_n, & \quad \tau^n_j := \delta_{\bar{y}_j} \otimes \sigma^{l_n}, \\
\Omega^n_{ij} := S^n_i \times T^n_{j} \subset X_n \times Y_n, & \quad \omega^n_{ij} := \sigma^n_i \otimes \tau^n_j
\end{align*}
for every $i, j = 0, 1$, and let $p_n \colon X_n \times Y_n \to Z$ be the map defined by
\begin{equation*}
p_n(x) := z_{ij} \text{ if } x \in \Omega^n_{ij}.
\end{equation*}
Note that ${p_n}_* (m_{X_n} \otimes m_{Y_n}) = m_Z$. Let $\varepsilon > 0$ be a sufficiently small arbitrary real number. We find a number $N \in \N$ such that 
\begin{align*}
& |F_n(\varepsilon, \varepsilon) - F(\varepsilon, \varepsilon)| < \varepsilon,  \\
& |\alpha_n - \alpha| < \varepsilon, \quad |\beta_n - \beta| < \varepsilon, \quad |\gamma_n - \gamma| < \varepsilon, \\
& \OD((S^n_i, \| \cdot \|, \sigma^n_i); -\varepsilon) < \varepsilon, \\
& \text{and } \OD((T^n_j, \| \cdot \|, \tau^n_i); -\varepsilon) < \varepsilon 
\end{align*}
hold for any $n \geq N$ and $i, j = 0, 1$. Let us prove that $p_n$ enforces $(72F(\varepsilon, \varepsilon) + 75\varepsilon)$-concentration of $X_n \times_{F_n} Y_n$ to $Z$ for any $n \geq N$. Fix $n$ with $n \geq N$. We first prove that $p_n$ is 1-Lipschitz up to $\varepsilon$. Take any $x_0, x'_0 \in S^n_0$, $x_1, x'_1 \in S^n_1$, $y_0, y'_0 \in T^n_0$, and $y_1, y'_1 \in T^n_1$ . We see that
\begin{align*}
s \leq \|x_0 - x_1\| \leq s_n, \quad 0 \leq \|x_i - x'_i\| \leq 2r_n, \\
t \leq \|y_0 - y_1\| \leq t_n, \quad 0 \leq \|y_j - y'_j\| \leq 2\rho_n.
\end{align*}
Thus we have
\begin{equation*}
d_{F_n}((x_0, y_j), (x_1, y'_j)) = F_n(\|x_0 - x_1\|, \|y_j - y'_j \|) \geq \alpha_n \geq \alpha -\varepsilon = d_Z(z_{0j}, z_{1j}) -\varepsilon
\end{equation*}
for every $j = 0, 1$, and
\begin{equation*}
d_{F_n}((x_i, y_0), (x'_i, y_1)) = F_n(\|x_i - x'_i\|, \|y_0 - y_1\|) \geq \beta_n \geq  \beta - \varepsilon = d_Z(z_{i0}, z_{i1}) - \varepsilon
\end{equation*}
for every $i = 0, 1$, and
\begin{equation*}
d_{F_n}((x_0, y_j), (x_1, y_{1-j})) = F_n(\|x_0 - x_1\|, \|y_0 - y_1\|) \geq \gamma_n \geq \gamma - \varepsilon = d_Z(z_{0j}, z_{1,1-j}) - \varepsilon
\end{equation*}
for every $j = 0, 1$. These imply that the map $p_n$ is 1-Lipschitz up to $\varepsilon$. By Lemma \ref{mmg5.27}, we have
\begin{equation*}
{p_n}^* \Lip_1(Z) \subset U_{2\varepsilon}(\Lip_1(X_n \times_{F_n} Y_n)).
\end{equation*}
We prove the other side inclusion. We take any function $f_n \in \Lip_1(X_n \times_{F_n} Y_n)$ and define a function $g_n : Z \to \R$ by
\begin{equation*}
g_n(z_{ij}) := \lm(f_n; \omega^n_{ij})
\end{equation*}
for every $i, j = 0, 1$. By Lemma \ref{mmg7.31} and Lemma \ref{key_1dim},
\begin{align*}
& \LR((\Omega^n_{ij}, d_{F_n}, \omega^n_{ij}); -2\varepsilon) \\
\leq & \OD((\Omega^n_{ij}, d_{F_n}, \omega^n_{ij}); -2\varepsilon) \\
\leq & 4F_n(\OD((S^n_i, \| \cdot \|, \sigma^n_i); -\varepsilon), 0) \\
& + 8F_n(0, \OD((T^n_j, \| \cdot \|, \tau^n_j); -\varepsilon))\\
\leq & 8F_n(\varepsilon, \varepsilon) + 16F_n(\varepsilon, \varepsilon) = 24F_n(\varepsilon, \varepsilon).
\end{align*}
Thus we have 
\begin{align*}
& m_{X_n} \otimes m_{Y_n} \left\{z \in X_n \times Y_n \midd |f_n(z) - {p_n}^* g_n(z)| > 24F_n(\varepsilon, \varepsilon) \right\} \\
& = \frac{1}{4} \sum_{i, j=0}^1  \omega^n_{ij} \left\{z \in \Omega^n_{ij} \midd |f_n(z) - \lm(f_n; \omega_{ij})| > 24F_n(\varepsilon, \varepsilon) \right\} \leq 2\varepsilon,
\end{align*}
which implies $\kf^{m_{X_n}\otimes m_{Y_n}}(f_n, {p_n}^*g_n) < 24F_n(\varepsilon, \varepsilon) + 2\varepsilon < 24F(\varepsilon, \varepsilon) + 26\varepsilon$. 

Let us prove the following three inequalities. For every $i, j = 0, 1$,
\begin{align}
|\lm(f_n; \omega^n_{0j}) - \lm(f_n; \omega^n_{1j})| & < \alpha + 48F(\varepsilon, \varepsilon) + 49\varepsilon, \label{g_alpha}\\
|\lm(f_n; \omega^n_{i0}) - \lm(f_n; \omega^n_{i1})| & < \beta + 48F(\varepsilon, \varepsilon) + 49\varepsilon, \label{g_beta}\\
|\lm(f_n; \omega^n_{0j}) - \lm(f_n; \omega^n_{1,1- j})| & < \gamma + 48F(\varepsilon, \varepsilon) + 49\varepsilon. \label{g_gamma}
\end{align}
We define angles $\theta_n, \theta'_n, \theta''_n, \phi_n, \phi'_n, \phi''_n \in [0, \pi]$ by
\begin{equation*}
\left\{\begin{array}{l} \cos{\theta_n} := \frac{1}{2r_n^2}\left( \frac{s_n^2 + s^2}{2} - u_{\alpha_n}^2\right), \\ \cos{\theta'_n} := \frac{1}{2r_n^2}\left( 2r_n^2 - u_{\beta_n}^2\right), \\ \cos{\theta''_n} := \frac{1}{2r_n^2}\left( \frac{s_n^2 + s^2}{2} - u_{\gamma_n}^2\right), \end{array}\right. \left\{\begin{array}{l} \cos{\phi_n}  := \frac{1}{2\rho_n^2}\left( 2\rho_n^2 - v_{\alpha_n}^2\right),  \\  \cos{\phi'_n}  := \frac{1}{2\rho_n^2}\left( \frac{t_n^2 + t^2}{2} - v_{\beta_n}^2\right), \\ \cos{\phi''_n}  := \frac{1}{2\rho_n^2}\left( \frac{t_n^2 + t^2}{2} - v_{\gamma_n}^2\right), \end{array}\right.
\end{equation*}
where $(u_{\alpha_n}, v_{\alpha_n})$, $(u_{\beta_n}, v_{\beta_n})$, and $(u_{\gamma_n}, v_{\gamma_n})$ are minimizers of $\alpha_n$, $\beta_n$ and $\gamma_n$ respectively, that is, $u_{\alpha_n}, u_{\gamma_n} \in [s, s_n]$, $u_{\beta_n} \in [0, 2r_n]$, $v_{\beta_n}, v_{\gamma_n} \in [t, t_n]$, $v_{\alpha_n} \in [0, 2\rho_n]$ such that 
\begin{equation*}
F_n(u_{\alpha_n}, v_{\alpha_n}) = \alpha_n, \ F_n(u_{\beta_n}, v_{\beta_n}) = \beta_n, \text{ and } F_n(u_{\gamma_n}, v_{\gamma_n}) = \gamma_n.
\end{equation*}
In addition, we define rotations $\Theta_n, \Theta'_n, \Theta''_n \colon \R^{k_n+1} \to \R^{k_n+1}$ by
\begin{align*}
\Theta_n^*&(a_1, \ldots, a_{k_n+1}) \\
:= & (a_1 \cos{\theta_n^*} - a_2 \sin{\theta_n^*}, a_1 \sin{\theta_n^*} + a_2 \cos{\theta_n^*}, \ldots,  \\
& \ a_{k_n} \cos{\theta_n^*} - a_{k_n+1} \sin{\theta_n^*}, a_{k_n} \sin{\theta_n^*} + a_{k_n+1} \cos{\theta_n^*})
\end{align*}
for $(a_1, \ldots, a_{k_n+1}) \in \R^{k_n+1}$, where $(\Theta_n^*, \theta_n^*)$ is each of $(\Theta_n, \theta_n)$, $(\Theta'_n, \theta'_n)$, and $(\Theta''_n, \theta''_n)$, and we recall that $k_n$ is odd. We also define rotations $\Phi_n, \Phi'_n, \Phi''_n \colon \R^{l_n+1} \to \R^{l_n+1}$ by the same way as above $\Theta_n^*$, that is,
\begin{align*}
\Phi_n^*&(b_1, \ldots, b_{l_n+1}) \\
:= & (b_1 \cos{\phi_n^*} - b_2 \sin{\phi_n^*}, b_1 \sin{\phi_n^*} + b_2 \cos{\phi_n^*}, \ldots,  \\
& \ b_{l_n} \cos{\phi_n^*} - b_{l_n+1} \sin{\phi_n^*}, b_{l_n} \sin{\phi_n^*} + b_{l_n+1} \cos{\phi_n^*}).
\end{align*}
We first prove \eqref{g_alpha}. Fix $j \in \{0, 1\}$. Let $T_n \colon \Omega^n_{0j} \to \Omega^n_{1j}$ be the map defined by
\begin{equation*}
T_n\left((\bar{x}_0, a), (\bar{y}_j, b)\right) := \left((\bar{x}_1, \Theta_n(a)), (\bar{y}_j, \Phi_n(b))\right)
\end{equation*}
for $a \in S^{k_n}(r_n)$, $b \in S^{l_n}(\rho_n)$. Note that ${T_n}_* \omega^n_{0j} = \omega^n_{1j}$. For any $a \in S^{k_n}(r_n)$, $b \in S^{l_n}(\rho_n)$, we have
\begin{align*}
&\|(\bar{x}_0, a) - (\bar{x}_1, \Theta_n(a)) \| = \sqrt{s^2 + 2r_n^2(1 - \cos{\theta_n})} = u_{\alpha_n}, \\
&\|(\bar{y}_j, b) - (\bar{y}_j, \Phi_n(b))\| = \sqrt{2\rho_n^2(1 - \cos{\phi_n})} = v_{\alpha_n},
\end{align*}
so that for any $z \in \Omega^n_{0j}$,
\begin{equation*}
d_{F_n}(z, T_n(z)) = F_n(u_{\alpha}, v_{\alpha}) = \alpha_n.
\end{equation*}
The measure $(\id, T_n)_* (\omega^n_{0j})$ is an $\alpha_n$-transport plan between $\omega^n_{0j}$ and $\omega^n_{1j}$. By Lemma \ref{lm_lem} and Lemma \ref{key_1dim}, we have
\begin{align*}
& |\lm(f_n; \omega^n_{0j}) - \lm(f_n; \omega^n_{1j})|  \\
& \leq \alpha_n + \sum_{i = 0}^1 \OD((\Omega^n_{ij}, d_{F_n}, \omega^n_{ij}); -2\varepsilon) \leq \alpha_n + \sum_{i = 0}^1 24F_n(\varepsilon, \varepsilon) \\
& < \alpha + 48F(\varepsilon, \varepsilon) + 49\varepsilon.
\end{align*}
Thus \eqref{g_alpha} is obtained. We next prove \eqref{g_beta}. We fix $i \in \{0, 1\}$, and define a map $T'_n \colon \Omega^n_{i0} \to \Omega^n_{i1}$ by
\begin{equation*}
T_n\left((\bar{x}_i, a), (\bar{y}_0, b)\right) := \left((\bar{x}_i, \Theta'_n(a)), (\bar{y}_1, \Phi'_n(b))\right)
\end{equation*}
for $a \in S^{k_n}(r_n)$, $b \in S^{l_n}(\rho_n)$. Note that ${T'_n}_* \omega^n_{i0} = \omega^n_{i1}$. For any $a \in S^{k_n}(r_n)$, $b \in S^{l_n}(\rho_n)$, we have
\begin{align*}
&\|(\bar{x}_i, a) - (\bar{x}_i, \Theta'_n(a)) \| = \sqrt{2r_n^2(1 - \cos{\theta'_n})} = u_{\beta_n}, \\
&\|(\bar{y}_0, b) - (\bar{y}_1, \Phi'_n(b))\| = \sqrt{t^2 + 2\rho_n^2(1 - \cos{\phi'_n})} = v_{\beta_n},
\end{align*}
so that for any $z \in \Omega^n_{i0}$,
\begin{equation*}
d_{F_n}(z, T'_n(z)) = F_n(u_{\beta_n}, v_{\beta_n}) = \beta_n.
\end{equation*}
The measure $(\id, T'_n)_* (\omega^n_{i0})$ is a $\beta_n$-transport plan between $\omega^n_{i0}$ and $\omega^n_{i1}$. By Lemma \ref{lm_lem} and Lemma \ref{key_1dim}, we have
\begin{align*}
& |\lm(f_n; \omega^n_{i0}) - \lm(f_n; \omega^n_{i1})|  \\
& \leq \beta_n + \sum_{j = 0}^1 \OD((\Omega^n_{ij}, d_{F_n}, \omega^n_{ij}); -2\varepsilon) \\
& < \beta + 48F(\varepsilon, \varepsilon) + 49\varepsilon.
\end{align*}
We obtain \eqref{g_beta}. Let us next prove \eqref{g_gamma}. Fix $j \in \{0, 1\}$ and let $T''_n \colon \Omega^n_{0j} \to \Omega^n_{1, 1-j}$ be the map defined by
\begin{equation*}
T''_n\left((\bar{x}_0, a), (\bar{y}_j, b)\right) := \left((\bar{x}_1, \Theta''_n(a)), (\bar{y}_{1-j}, \Phi''_n(b))\right)
\end{equation*}
for $a \in S^{k_n}(r_n)$, $b \in S^{l_n}(\rho_n)$. Note that ${T''_n}_* \omega^n_{0j} = \omega^n_{1,1-j}$. For any $a \in S^{k_n}(r_n)$, $b \in S^{l_n}(\rho_n)$, we have
\begin{align*}
&\|(\bar{x}_0, a) - (\bar{x}_1, \Theta''_n(a)) \| = \sqrt{s^2 + 2r_n^2(1 - \cos{\theta''_n})} = u_{\gamma_n}, \\
&\|(\bar{y}_j, b) - (\bar{y}_{1-j}, \Phi''_n(b))\| = \sqrt{t^2 + 2\rho_n^2(1 - \cos{\phi''_n})} = v_{\gamma_n},
\end{align*}
so that for any $z \in \Omega^n_{0j}$,
\begin{equation*}
d_{F_n}(z, T''_n(z)) = F_n(u_{\gamma_n}, v_{\gamma_n}) = \gamma_n.
\end{equation*}
The measure $(\id, T''_n)_* (\omega^n_{0j})$ is a $\gamma_n$-transport plan between $\omega^n_{0j}$ and $\omega^n_{1,1-j}$. By Lemma \ref{lm_lem} and Lemma \ref{key_1dim}, we have
\begin{equation*}
|\lm(f_n; \omega^n_{0j}) - \lm(f_n; \omega^n_{1,1-j})| < \gamma + 48F(\varepsilon, \varepsilon) + 49\varepsilon.
\end{equation*}
We now complete the proof of three inequalities \eqref{g_alpha} -- \eqref{g_gamma}. Combining  \eqref{g_alpha} -- \eqref{g_gamma} and Lemma \ref{MHext} implies $\kf^{m_Z}(g_n, \Lip_1(Z)) < 48F(\varepsilon, \varepsilon) + 49\varepsilon$. Thus we have
\begin{equation*}
\begin{split}
& \kf^{m_{X_n} \otimes m_{Y_n} }(f_n, {p_n}^*\Lip_1(Z)) \\
& \leq \kf^{m_{X_n} \otimes m_{Y_n}}(f_n, {p_n}^*g_n) + \kf^{m_Z}(g_n, \Lip_1(Z)) < 72F(\varepsilon, \varepsilon) + 75\varepsilon.
\end{split}
\end{equation*}
Therefore the map $p_n$ enforces $(72F(\varepsilon, \varepsilon) + 75\varepsilon)$-concentration of $X_n \times_{F_n} Y_n$ to $Z$ for every $n \geq N$. By Theorem \ref{equiconc}, the sequence $\{X_n \times_{F_n} Y_n\}_{n \in \N}$ concentrates to $Z$ as $n \to \infty$. The proof of the claim is now completed.
\end{proof}
Since
\begin{equation*}
\alpha \leq F(s, 0), \ \beta \leq F(0, t), \text{ and } \gamma \leq F(s, t) - \eta,
\end{equation*}
the mm-space $X \times_F Y$ is not mm-isomorphic to $Z$. Thus Claim \ref{claim_nece} means that the condition (\ref{main1}) of Theorem \ref{main} does not hold. Therefore we obtain the implication from (\ref{main1}) to (\ref{main2}) of Theorem \ref{main}.
\end{proof}

\begin{rem}
In the above proof of Theorem \ref{main}, if there exist finite limits of both $\{s_n\}_{n \in \N}$ and $\{t_n\}_{n \in \N}$, then the three sequences $\{\alpha_n\}_{n \in \N}$, $\{\beta_n\}_{n \in \N}$, and $\{\gamma_n\}_{n \in \N}$ converge without taking a subsequence and these limits are
\begin{equation*}
\alpha = \min_{\substack{s \leq u_1 \leq s_\infty \\ 0 \leq v_1 \leq 2\rho}} F(u_1, v_1), \quad \beta = \min_{\substack{0 \leq u_2 \leq 2r \\ t \leq v_2 \leq t_\infty}} F(u_2, v_2), \quad \gamma = \min_{\substack{s \leq u_3 \leq s_\infty \\ t \leq v_3 \leq t_\infty}} F(u_3, v_3),
\end{equation*}
where 
\begin{align*}
s_\infty := \lim_{n\to\infty} s_n, \ t_\infty := \lim_{n\to\infty} t_n, \ r := \frac{\sqrt{s_\infty^2 - s^2}}{2}, \ \rho := \frac{\sqrt{t_\infty^2 - t^2}}{2}.
\end{align*}
\end{rem}

\section{Product of $N$ metric measure spaces}

In this section, we consider the concentration of product spaces of $N$ mm-spaces. Indeed, we generalize Theorem \ref{main} to the following.

\begin{thm}\label{main_Ndim}
Let $F_n, F \in \F^N$, $n = 1, 2, \ldots$. Assume that $F_n$ converges pointwise to $F$ as $n \to \infty$. Then the following conditions are equivalent to each other.
\begin{enumerate}
\item For any $N$ sequences $\{X_n^i\}_{n \in \N}$, $i = 1, \ldots, N$, of mm-spaces concentrating to mm-spaces $X^i$ respectively, the sequence $\{(\prod_{i=1}^N X_n^i, \\ d_{F_n}, \otimes_{i = 1}^N m_{X_n^i} )\}_{n \in \N}$ of their product spaces concentrates to the product space $(\prod_{i=1}^N X^i, d_F, \otimes_{i = 1}^N m_{X^i} )$ as $n \to \infty$.
\item For any $(s_1, \ldots, s_N) \in [0, +\infty)^N$, 
\begin{equation*}
\lim_{n \to \infty} (F_n(s_1, \ldots, s_N) - \inf_{s_i \leq s'_i} F_n(s'_1, \ldots, s'_N)) = 0.
\end{equation*}
\item For any $D > 0$, 
\begin{equation*}
\lim_{n \to \infty} \sup_{0 \leq s_i \leq D} (F_n(s_1, \ldots, s_N) - \inf_{s_i \leq s'_i} F_n(s'_1, \ldots, s'_N)) = 0.
\end{equation*}
\end{enumerate}
\end{thm}

We denote by $(\prod_{i=1}^N X^i)_p$ the $l_p$-product spaces of $X^1, \ldots, X^N$, which are generated by
\begin{equation*}
F_p^N(s_1, \ldots, s_N) := \left\{ \begin{array}{ll} \left(\sum_{i=1}^N s_i^p\right)^{\frac{1}{p}} & \text{if } p < +\infty, \\  \displaystyle \max_{i=1, \ldots, N} s_i & \text{if } p = +\infty.  \end{array}\right.
\end{equation*}

\begin{lem}\label{key_lp_N}
Let $p \in [1, +\infty]$ and let $X^1, \ldots, X^N$ be $N$ mm-spaces. Then we have
\begin{equation}\label{key_lp_N_eq}
\begin{split}
& \OD((\prod_{i=1}^N X^i)_p; - \sum_{i=1}^N\kappa_i) \\
& \leq \OD(X^1 ; - \kappa_1) + 2 \sum_{i=2}^N \OD(X^i ; -\kappa_i)
\end{split}
\end{equation}
for any $\kappa_1 \in (0, 1)$ and any $\kappa_2, \ldots, \kappa_N \in (0, 1/2)$.
\end{lem}

\begin{lem}\label{key_F_N}
Let $F \in \F^N$ and let $X^1, \ldots, X^N$ be $N$ mm-spaces. Then
\begin{equation}\label{key_F_N_eq}
\begin{split}
& \OD((\prod_{i=1}^N X^i, d_F, \bigotimes_{i = 1}^N m_{X^i} ) ; - 2 \sum_{i=1}^N\kappa_i) \\
& \leq 4F^1 (\OD(X^1 ; - \kappa_1)) + 8\sum_{i=2}^N F^i (\OD(X^i ; -\kappa_i))
\end{split}
\end{equation}
for any $\kappa_1 \in (0, 1)$ and any $\kappa_2, \ldots, \kappa_N \in (0, 1/4)$, where $F^i := F \circ \iota_i$ and $\iota_i \colon [0, +\infty) \to [0, +\infty)^N$ is the natural $i$-th inclusion map.
\end{lem}

\begin{proof}[Proof of Lemma \ref{key_lp_N}]
Since the $l_p$-product has the iterated property
\begin{equation*}
(\prod_{i=1}^N X^i)_p = (\prod_{i=1}^{N-1} X^i)_p \times_p X^N,
\end{equation*}
by Lemma \ref{key_lp}, we have
\begin{equation*}
\begin{split}
& \OD((\prod_{i=1}^N X^i)_p; - \sum_{i=1}^N\kappa_i) \\
& \leq \OD((\prod_{i=1}^{N-1} X^i)_p ; - \sum_{i=1}^{N-1}\kappa_i) + 2 \OD(X^N ; -\kappa_N).
\end{split}
\end{equation*}
We obtain \eqref{key_lp_N_eq} by the induction.
\end{proof}

\begin{proof}[Proof of Lemma \ref{key_F_N}]
In the same way as the proof of Lemma \ref{key_F}, since $F \leq \sum_{i=1}^N F^i := G$ and the mm-space $(\prod_{i=1}^N X^i, d_G, \otimes_{i = 1}^N m_{X^i})$ is mm-isomorphic to
\begin{equation*}
(\prod_{i=1}^N (X^i, F^i \circ d_{X^i}, m_{X^i}) )_1,
\end{equation*}
by Lemma \ref{key_lp_N}, we have \eqref{key_F_N_eq}.
\end{proof}

\begin{proof}[Proof of Theorem \ref{main_Ndim}]
We are able to prove Theorem \ref{main_Ndim} by imitating the proof of Theorem \ref{main}. We obtain the implication from (2) to (3) by imitating Lemmas \ref{ptcpt} and \ref{equiv_main}. In the proof of the implication from (3) to (1), the key to imitate is an estimate of the observable diameter of the product space. We have already obtained Lemma \ref{key_F_N}, so that we obtain this implication. Let us prove the implication from (1) to (2). Assume that the condition (2) does not hold. Up to choosing a subsequence of $n$, we are able to assume that there exist a real numbers $\eta > 0$, an $N$-tuple $(s^1, \ldots, s^N) \in [0, +\infty)$ and a sequence $\{(s_n^1, \ldots, s_n^N)\}_{n \in \N} \subset [0, +\infty)$ such that 
\begin{equation*}
s^i < s_n^i \text{ and } F_n(s^1, \ldots, s^N) > F_n(s_n^1, \ldots, s_n^N) + \eta
\end{equation*}
for any $n \in \N$ and every $i$. Let $X^{i}$, $i = 1, \ldots, N$, be the $N$ mm-spaces defined by
\begin{equation*}
X^i :=(\{x_0^i, x_1^i\}, d_{X^i}, \frac{1}{2}\delta_{x_0^i}+\frac{1}{2}\delta_{x_1^i}), \quad d_{X^i}(x_0^i, x_1^i) := s^i.
\end{equation*}
We consider the $N$ mm-spaces $X_n^i$ defined by
\begin{equation*}
X_n^i := X^{i} \times_2 (S^{k_n^i}(r_n^i), \|\cdot\|, \sigma^{k_n^i}),
\end{equation*}
where $r_n^i := \sqrt{(s_n^{i})^2 - (s^{i})^2}/2$ and $k_n^i := 2\max\{n, \lceil (r_n^i)^4\rceil\} + 1$. By imitating the proof of Claim \ref{claim_nece}, we see that $\{(\prod_{i=1}^N X_n^i, d_{F_n}, \otimes_{i=1}^N m_{X_n^i})\}_{n \in \N}$ has a subsequence that does not concentrate to $(\prod_{i=1}^N X^i, d_F, \otimes_{i=1}^N m_{X^i})$. Thus we obtain the implication from (1) to (2). The proof is completed.
\end{proof}

\begin{ex}
The following function is an example of metric preserving functions that are not the iterated type.
\begin{equation*}
F_{\text{cyc}}(s_1, s_2, s_3) := \max\{s_1 + s_2, s_2 + s_3, s_3 + s_1\}.
\end{equation*}
$F_{\text{cyc}}$ does not have the iterated property like $F_p^N$. Theorem \ref{main_Ndim} can be applied to such functions.
\end{ex}

\end{document}